\documentclass{amsart}
\usepackage{mystyle}
\usepackage[utf8]{inputenc}
\usepackage{bm}
\usepackage{accents}

\newcommand{\myoset}[2]{\accentset{\scriptstyle {#1}}{#2}}

\title{The Dupont Homotopy Formula and Stellar Subdivision}
\author{Benjamin I. Albert}

\begin{document}
\begin{abstract}
  The Dupont homotopy, a classical construction in the algebraic topology of triangulated smooth manifolds, has been revived in the last decade in the construction of an effective field theory where it appears as a propagator.  In this paper, we ask a question of relevance to the renormalization group of this theory: is Dupont's construction compatible with stellar subdivision?
\end{abstract}

\maketitle
\tableofcontents

\section{Introduction}
\label{sec:introduction}
Whitney realized \cite{Whitney_1957} that for any smooth manifold $M$ with triangulation $|\Xi| \cong M$, there is a cochain map $W: C^{\bullet}(\Xi) \to \Omega^{\bullet}_{\Xi}(M)$ taking simplicial cochains to piecewise smooth differential forms, with $W$ being a section for (i.e.\ a right inverse to) the integration map $R: \Omega_{\Xi}^{\bullet}(M) \to C^{\bullet}(M)$.  Because, by the de Rham isomorphism theorem, $R$ induces an isomorphism on cohomology, $W$ also induces an isomorphism on cohomology, so the image of $W$, the space of Whitney forms, generates the de Rham cohomology of $M$.

Several decades later, Dupont \cite{Dupont_1978} proved the stronger result that there is a deformation retraction of $\Omega_{\Xi}^{\bullet}(M)$ onto $C^{\bullet}(\Xi)$.  That is, he showed there exists a homotopy $s$ between the identity map $1$ and $WR$.  Dupont was interested in the study of characteristic classes and constructing a universal Chern-Weil homomorphism taking an invariant polynomial on the Lie algebra of a Lie group $G$ to a cohomology class of the classifying space $|BG|$.  Dupont used his homotopy, in particular, as an intermediate step to relate the de Rham complex of $BG$, a simplicial manifold, to the simplicial cochain complex of its geometric realization $|BG|$.  Dupont was then able to relate the universal characteristic class he constructed as an element of the de Rham cohomology of $BG$ to the usual notion of universal characteristic class as an element of the simplicial cohomology of $|BG|$.

More recently, Getzler \cite{2004math......4003G} has made use of the Dupont homotopy in his study of nilpotent Lie algebras and more generally nilpotent $L_{\infty}$-algebras.  He used the Dupont homotopy to construct a space $\gamma(\mathfrak{g})$ that is isomorphic to $BG$ when $\mathfrak{g}$ is a nilpotent Lie algebra.  Since $\gamma(\mathfrak{g})$ exists for any nilpotent $L_{\infty}$-algebra $\mathfrak{g}$, it can be thought of as a generalized notion of classifying space.

In a different direction, Mnev \cite{2008arXiv0809.1160M} used the Dupont homotopy as a propagator for BF theory on triangulated manifolds.  In his paper, the effective action, a functional on the space of Lie-algebra-valued simplicial cochains, is calculated for a variety of familiar topological spaces, and explicit combinatorial formulas are written down.  The paper also begins the study of how to glue Dupont homotopies when gluing triangulated manifolds with boundary.  In a subsequent paper, Cattaneo, Mnev, and Reshetikhin \cite{2017arXiv170105874C} treat the more general setting of a cellular complex where one is required to make a noncanonical choice of a deformation retraction of a suitable space of differential forms onto the space of cellular cochains.  They construct the effective BF action on cellular manifolds with boundary and show that the construction is compatible with gluing of cobordisms.

By a theorem of Alexander (see Lickorish \cite{1999math.....11256L} for a modern proof), any two triangulations of a topological manifold (with a common refinement) are related by a sequence of stellar subdivisions and stellar weldings (i.e.\ inverse stellar subdivisions).  To motivate the results of our paper, we make a few remarks on how stellar subdivision solves an elementary problem in algebraic topology, which is to show that the simplicial chains of any two triangulations are homotopy equivalent.  As is standard practice, this problem can be solved by appealing to the theory of singular chains (see Hatcher's book \cite{Hatcher_2002}).  However, if we require the triangulations have a common refinement, there is a more direct approach.  In this case, it suffices to be able to exhibit for any triangulation $\Xi$, a deformation retraction from the simplicial chains on a stellar subdivision $\star \Xi$ onto the simplicial chains of $\Xi$.  From a sequence of stellar subdivisions and stellar weldings relating the two triangulations, we then would have a sequence of homotopy equivalences of simplicial chain complexes which can be composed.  

Having constructed the deformation retraction of $C_{\bullet}(\star \Xi)$ onto $C_{\bullet}(\Xi)$ of the previous paragraph, by dualizing the retraction map, we have the ability to include $C^{\bullet}(\Xi)$ as a subspace of $C^{\bullet}(\star\Xi)$, and thereby, under the Whitney map, also define a subspace of the Whitney forms on $M$ for the triangulation $\star\Xi$.  A basic question which arises is whether this subspace is the same as the space of Whitney forms on $M$ for the triangulation $\Xi$.  Then the more involved natural question to ask is whether the Whitney map, the integration map and the Dupont homotopy are compatible in some appropriate way with stellar subdivision of the triangulation.  

Further motivation for our compatibility question comes from quantum field theory on a lattice, where the renormalization group should relate the physics on a lattice to the physics on a refinement of that lattice.  By ``integrating out'' the additional degrees of freedom of the fields on the refined lattice, one recovers the physics depending on the fields on the original lattice.  This is the renormalization group picture pioneered by Kadanoff and Wilson (see \cite{Kopietz_2010}, Chapter 3).  In the setting of the simplicial BF theory of Mnev \cite{2008arXiv0809.1160M} the renormalization group should be interpreted as a compatibility between the effective actions associated to different triangulations (with common refinement), or, as suffices by Alexander's theorem, just a triangulation $\Xi$ and any stellar subdivision $\star \Xi$.  This compatibility condition might be a statement about propagators: The sum of the propagator used to construct the effective action for $\star \Xi$ (the Dupont homotopy for $\star \Xi$) and the propagator used to integrate out the additional degrees of freedom for fields on $\star \Xi$ (the stellar welding homotopy) is equal to the propagator used to construct the effective action for $\Xi$ (the Dupont homotopy for $\Xi$).  The main theorem in this paper, Theorem \ref{thm:dup_stell_subdiv}, gives a precise but slightly weaker version of the above (naive) statement.

As a broad outline of the paper, in Section \ref{sec:dup_hom}, we review the construction of the integration map, the Whitney map, and the Dupont homotopy and show that these maps together specify a deformation retraction of $\Omega^{\bullet}_{\Xi}(M)$ onto $C^{\bullet}(\Xi)$.  The exposition essentially follows the classic monograph by Dupont \cite{Dupont_1978} and the paper of Getzler \cite{2004math......4003G}.  In Section \ref{sec:cub_dup_hom}, we construct the Dupont homotopy for cubical chains starting from the $1$-dimensional Dupont homotopy following the tensor product construction for deformation retractions.

In Section \ref{sec:stellar_subdiv_form}, we define a deformation retraction from simplicial (co)chains on $\star\Delta^n$, a stellar subdivision of the $n$-simplex, onto simplicial (co)chains on $\Delta^n$, the $n$-simplex.  From this, we induce a deformation retraction from simplicial (co)chains on a simplicial complex $\star\Xi$ onto simplicial (co)chains on its stellar welding $\Xi$.  In Section \ref{sec:cub_stell_subdiv}, we define the notion of cubical subdivision.  Employing the tensor product construction, we find a deformation retraction from cubical (co)chains on $\star\square^n$, a cubical subdivision of the $n$-cube, onto the cubical (co)chains on the $n$-cube $\square^n$.  From this, we induce a deformation retraction from cubical (co)chains on a cubical complex $\star\Xi$ onto cubical cochains on its cubical welding $\Xi$.

Section \ref{sec:comp-stell-dup} lays out the main compatibility results of the paper.  The deformation retraction from $\Omega^{\bullet}_{\star \Xi}(M)$ onto $C^{\bullet}(\star \Xi)$ can be composed with the stellar welding deformation retraction from $C^{\bullet}(\star \Xi)$ onto $C^{\bullet}(\Xi)$.  We find that the restriction of this composed deformation retraction is equivalent in a certain sense to the deformation retraction of $\Omega^{\bullet}_{\Xi}(M)$ onto $C^{\bullet}(\Xi)$.  In Section \ref{sec:comp_cub_stell}, we examine the cubical case of this compatibility result.

Lastly, in Sections \ref{sec:el_coll} and \ref{sec:stell_subdiv2}, we recover the stellar welding deformation retraction whose formula had been posited in Section \ref{sec:stellar_subdiv_form} without motivation.  On a historical note, this calculation was one of the main things that inspired us to work on the construction of Section \ref{sec:comp-stell-dup}.  In Section \ref{sec:el_coll} we recall the elementary collapse deformation retraction (whose formula can be found in \cite{2017arXiv170105874C}).  In Section \ref{sec:stell_subdiv2}, we show that a stellar subdivision can be constructed by a sequence of elementary expansions (inverse elementary collapses) followed by elementary collapses.  For stellar subdivision at a $k$-simplex, we exhibit $k + 1$ such sequences.  Each sequence gives rise to a zigzag of elementary collapse deformation retractions, whose composition is still, in fact, a deformation retraction.  We prove that the average of these $k + 1$ deformation retractions is equal to our stellar welding deformation retraction.

\section{Dupont Homotopy Formula}
\label{sec:dup_hom}

\subsection{Definitions}
\label{sec:defs}
Let $\Delta^n$ denote the standard $n$-simplex, thought of as a simplicial complex whose geometric realization is given by
\begin{align*}
  |\Delta^n| = [e_0, \dots, e_n] =  \lbrace t_0e_0 + \dots + t_ne_n : \sum_it_i = 1 \rbrace \subset \RR^{n + 1}.
\end{align*}

For any simplicial chain $\alpha \in C_{\bullet}(\Delta^n)$, let $\widehat{\alpha} \in C^{\bullet}(\Delta^n)$ be its dual simplicial cochain.
\begin{mydef}
  For each simplex $[i_0, \dots, i_p]$ in $\Delta^n$, define the \emph{Whitney form} $\overline{\omega}_{i_0, \dots, i_p}$, a $p$-form on $|\Delta^n|$, by $\overline{\omega}_{i_0, \dots, i_p} = p! \, \omega_{i_0, \dots, i_p}$, where
  \begin{align*}
    \omega_{i_0, \dots, i_p} = \sum_l (-1)^l t_{i_l}dt_{i_0} \dots \widehat{dt_{i_l}} \dots dt_{i_n}.
  \end{align*}
  Then define the \emph{Whitney map} $W:C^{\bullet}(\Delta^n) \to \Omega^{\bullet}(|\Delta^n|)$ by
\begin{align*}
  W(\widehat{\alpha}) &= \overline{\omega}_{i_0, \dots, i_p}
\end{align*}
for $\alpha = [i_0, \dots, i_p]$.
\end{mydef}
\noindent Note that
\begin{align*}
  \sum_k \omega_{k, i_0, \dots, i_p} &= \left(\sum_k t_k\right) dt_{i_0} \dots dt_{i_p} -  \sum_l (-1)^l t_{i_l} \left(\sum_k dt_k\right) dt_{i_0} \dots \widehat{dt_{i_l}} \dots dt_{i_p}\\
  &= dt_{i_0} \dots dt_{i_p}
\end{align*}
Therefore $W$ is a cochain map:
\begin{align*}
  dW(\widehat{\alpha}) &= (p + 1)!\,dt_{i_{0}} \dots dt_{i_p} = \sum_k \overline{\omega}_{k, i_0, \dots, i_p} = W(d\widehat{\alpha})
\end{align*}
Recall that the integration map
  \begin{align*}
    R = \sum_{p = 0}^n\sum_{i_0 < \dots < i_p} \reallywidehat{[i_0, \dots, i_p]} \int_{[i_0, \dots, i_p]}
  \end{align*}
  is also a cochain map as a consequence of Stokes' Theorem.
  
  We shall verify that $\int_{[i_0, \dots, i_p]} \omega_{i_0, \dots, i_p} = \frac{1}{p!}$ in the next section.  This implies that
\begin{align*}
  RW(\reallywidehat{[i_0, \dots, i_p]}) &= \reallywidehat{[i_0, \dots, i_p]}\, p! \int_{[i_0, \dots, i_p]} \omega_{i_0, \dots, i_p} = \reallywidehat{[i_0, \dots, i_p]}.
\end{align*}
That is, $RW = 1$.  Dupont discovered that while $WR \neq 1$, there is a homotopy between $1$ and $WR$.

The Dupont homotopy is expressed in terms of the Whitney forms and certain maps $h^i$ of degree $-1$, for $i$ ranging from $0$ to $n$.  To define $h^i$, we first define the maps
\begin{align*}
  \phi^i: [0, 1] \times |\Delta^n| \to |\Delta^n|
\end{align*}
as the contraction of the $n$-simplex onto its $i$-th vertex $e_i$. 
Then $h^i= \pi_{*}(\phi^i)^{*}$ where $\pi_{*}$ is integration along the fiber $[0, 1]$.  

\begin{mydef}
  The \emph{Dupont homotopy} is given by the formula
  \begin{align*}
    s = -\sum_{k = 0}^{n - 1} \sum_{i_0 < \dots < i_k} \overline{\omega}_{i_0, \dots, i_k} h^{i_k}\dots h^{i_0} 
  \end{align*}
\end{mydef}
 
We now give some motivation for why this might be a homotopy.  At the beginning of Section \ref{sec:proofs} we shall argue that $h^i$ is a homotopy between the evaluation map at $e_i$, $\varepsilon^i: \Omega^{\bullet}(|\Delta^n|) \to \Omega^0(|\Delta^n|) \to \RR$, and the identity map $1: \Omega^{\bullet}(|\Delta^n|) \to \Omega^{\bullet}(|\Delta^n|)$.  In a different direction, in Section \ref{sec:comps}, 
where in particular we derive the basic properties of $h^i$, we will show that
\begin{align*}
  (-1)^p\varepsilon^{i_p}h^{i_{p - 1}}\dots h^{i_0}(\omega_{i_0, \dots, i_p}) = \frac{1}{p!}.
\end{align*}
We have also claimed that
\begin{align*}
  \int_{[i_0, \dots, i_p]}\omega_{i_0, \dots, i_p} = \frac{1}{p!}.
\end{align*}
This motivates the plausibility of a more general statement: for any $\omega \in \Omega^{\bullet}(\Delta^n)$,
  \begin{align*}
  \int_{[i_0, \dots, i_p]}\omega = (-1)^p\varepsilon^{i_p}h^{i_{p - 1}}\dots h^{i_0}(\omega).
  \end{align*}
We shall prove this as a lemma in Section \ref{sec:comps}.  

With $R$, $s$ and $W$ written in terms of $h^i$ and $\overline{\omega}_{i_0, \dots, i_k}$ and our knowledge of the properties of $h^i$ and $\omega_{i_0, \dots, i_k}$, the desired statement $ds + sd = 1 - WR$ seems plausible.  We shall prove this along with the properties $s^2 = 0$, $sW = 0$, and $Rs = 0$ (which make the deformation retraction ``special'') in Section \ref{sec:proofs}.  In summary:
\newtheorem*{thm:dupont}{Theorems \ref{thm:proof-main-theorem} \& \ref{thm:dupont-special-main-theorem}}
\begin{thm:dupont}
  There is a (special) deformation retraction of the differential forms $\Omega^{\bullet}(|\Delta^n|)$ on the $n$-simplex onto the simplicial cochains $C^{\bullet}(\Delta^n)$ on the $n$-simplex:
\begin{align*}
  \myDR{(C^{\bullet}(\Delta^n), d)}{(\Omega^{\bullet}(|\Delta^n|), d)}{W}{R}{s}
\end{align*}
where $W$ is the Whitney map, $R$ is the integration map, and $s$ is the Dupont homotopy. 
\end{thm:dupont}

The proof we shall give of Theorems \ref{thm:proof-main-theorem} \& \ref{thm:dupont-special-main-theorem} more or less follows the one in Getzler's paper \cite{2004math......4003G} which is mostly based the one outlined in Dupont's monograph \cite{Dupont_1978}.  The differences here are mainly greater overall detail.  Also, our approach to the proof of Lemma \ref{lem:integral} is perhaps more illustrative.  Lastly, we correct a few incorrect formulas in \cite{2004math......4003G} that result from that paper mostly following the conventions of \cite{Dupont_1978}, but changing the definition of $\phi^i$.

\subsection{Basic Properties}
Recall that, by definition, $\phi^i$ is a homotopy between the identity map on $\Delta^n$ and the constant map with value $e_i$; that is,
\begin{align*}
  \phi^i\left(s, \sum_jt_je_j\right) = (1 - s)\sum_jt_je_j + se_i = \sum_j ((1 - s) t_j + s\delta_{ij})e_j.
\end{align*}
\label{sec:comps}
  We compute in coordinates
\begin{align*}
  h^i(f dt_{i_{0}} \dots dt_{i_p}) &= \pi_{*} [f(\phi^i(s, \mathbf{t}))\, d((1 - s)t_{i_0}) \dots d((1 - s)t_i + s) \dots d((1 - s)t_{i_p})]\\
                                   &= \pi_{*} \left[ f(\phi^i(s, \mathbf{t}))(1 - s)^p \sum_j (-1)^j(\delta_{i, i_{j}} - t_{i_j})\,ds\, dt_{i_0} \dots  \widehat{dt_{i_j}} \dots dt_{i_p} \right]\\
  &= \left(\int_0^1(1 - s)^pf(\phi^i(s, \mathbf{t})) \, ds \right)  \sum_j (-1)^j(\delta_{i, i_{j}} - t_{i_j})\, dt_{i_0} \dots  \widehat{dt_{i_j}} \dots dt_{i_p}
\end{align*}
Introduce the vector fields $E_i$  on $|\Delta^n|$ given by
\begin{align*}
  E_i = \sum_j (\delta_{i, j} - t_j) \partial_{t_j}.
\end{align*}
A priori, these are only vector fields on $\RR^{n + 1}$.  However,
\begin{align*}
  E_i\left(\sum_kt_k - 1\right) = 1 - \sum_k t_k
\end{align*}
so $E_i$ preserves the ideal generated by $\sum_kt_k - 1$, which implies that $E_i$ is a vector field on $|\Delta^n|$.

The relation $(\phi^i_s)^{*}(\delta_{i,j} - t_j) = (1 - s)(\delta_{i,j} - t_j)$ implies that
\begin{align*}
  (\phi^i_s)^{*}\imath_{E_i}(f dt_{i_{0}} \dots dt_{i_p}) = f(\phi^i(s, \mathbf{t})) (1 - s)^{p + 1}\sum_j (-1)^j(t_{i_j} - \delta_{i, i_{j}})\, dt_{i_0} \dots  \widehat{dt_{i_j}} \dots dt_{i_p}
\end{align*}
and thus
\begin{myprop} \label{prop:hiomega_gen}
  For any form $\omega \in \Omega^{\bullet}(|\Delta^n|)$, we have
  \begin{align*}
    h^i(\omega) = \int_0^1 \frac{ds}{1 - s}(\phi^i_s)^{*}\imath_{E_i}\omega.
  \end{align*}
\end{myprop}
\noindent The next proposition, a useful computational tool, will imply that
\begin{align*}
(-1)^p\varepsilon^{i_p}h^{i_{p - 1}} \dots h^{i_0}\omega_{i_0, \dots, i_p} = 1/p!.
\end{align*}
\begin{myprop}\label{prop:hiomega}
  \begin{align*}
    \imath_{E_i}(\omega_{i_0, \dots, i_p}) &= (-1)^{l + 1}\omega_{i_0, \dots, \widehat{i_l}, \dots, i_p}.
  \end{align*}
  and therefore
  \begin{align*}
    h^i(\omega_{i_0, \dots, i_p}) = \frac{(-1)^{l + 1}}{p}\omega_{i_0, \dots, \widehat{i_l}, \dots, i_p}
  \end{align*}
  if $i = i_l$ for some $l$ and $h^i(\omega_{i_0, \dots, i_p}) = 0$ otherwise.
\end{myprop}
\begin{proof}
  Let $i_{\Delta^n}: |\Delta^n| \to \RR^{n + 1}$ denote the inclusion map.  From the definitions,
  \begin{align*}
    \imath_{E_i}(\omega_{i_0, \dots, i_p}) &=\imath_{E_i} i_{\Delta^n}^{*}\left(\sum_j t_{i_j} \imath_{\partial_{t_{i_j}}} dt_{i_0} \dots dt_{i_p}\right)\\
                                           &= i_{\Delta^n}^{*} \left(\sum_{j, k} t_{i_j}(\delta_{i, i_k} - t_{i_k}) \imath_{\partial_{t_{i_k}}}\imath_{\partial_{t_{i_j}}} (dt_{i_0} \dots dt_{i_p})\right)\\
                                           &= -i_{\Delta^n}^{*} \left(\sum_j t_{i_j}\imath_{\partial_{t_{i_j}}} \imath_{\partial_{t_i}}  (dt_{i_0} \dots dt_{i_p})\right)
  \end{align*}
which implies that
\begin{align*}
  (\phi^i_s)^{*}\imath_{E_i}(\omega_{i_0, \dots, i_p})  &= - (1 - s)^pi_{\Delta^n}^{*}\left(\sum_j t_{i_j}\imath_{\partial_{t_{i_j}}} \imath_{\partial_{t_i}}  (dt_{i_0} \dots dt_{i_p})\right)\\
  &=
    \begin{cases}
      (-1)^{l + 1}(1 - s)^p\omega_{i_0, \dots, \widehat{i_l}, \dots, i_p} & \text{ if $i = i_l$ for some $l$}\\
      0 & \text{ otherwise}.
    \end{cases}
\end{align*}
The result for $h^i(\omega_{i_0, \dots, i_p})$ now follows from Proposition \ref{prop:hiomega_gen}.
\end{proof}
Iterating the result of the above proposition, we find that
\begin{align*}
(-1)^p\varepsilon^{i_p}h^{i_{p - 1}} \dots h^{i_0}(\omega_{i_0, \dots, i_p}) = 1/p!.
\end{align*}
We would like to also show that $\int_{[i_0, \dots, i_p]} \omega_{i_0, \dots, i_p} = 1/p!$.  Rather than give a specific argument for this fact, we continue by arguing more generally that $\int_{[i_0, \dots, i_p]}\omega = (-1)^p\varepsilon^{i_p}h^{i_{p - 1}}\dots h^{i_0}(\omega)$ for any $p$-form $\omega$.

To calculate the integral, we pull back $\omega$ using a parametrization $[0, 1]^p \to [i_0, \dots, i_p] \subset |\Delta^n|$.  Let $\iota^{i_p}: [0, 1]^p \to [0,1]^p \times \Delta^n$ be the product of the identity map and the inclusion of the vertex $e_{i_p}$.  A natural candidate for the parametrization $[0, 1]^p \to [i_0, \dots, i_p]$ is
\begin{align*}
  F_{i_0, \dots, i_p} = \phi^{i_0} \circ (\id_{[0, 1]} \times \phi^{i_1}) \circ \dots \circ (\id_{[0, 1]^{p - 1}} \times \phi^{i_{p - 1}}) \circ \iota^{i_p}
\end{align*}
Note that for each $k \in \{0, \dots, p\}$ the image of $\phi^{i_k}|_{[0, 1] \times [i_{k + 1}, \dots, i_p]}$ is equal to $[i_k, \dots, i_p]$.  It follows that the image of $F_{i_0, \dots, i_p}$ is equal to $[i_0, \dots, i_p]$.  It turns out that $F_{i_0, \dots, i_p}$ is orientation preserving for $p$ even and orientation reversing for $p$ odd, which we prove below.  Let $R_p: [0, 1]^p \to [0, 1]^p$ be defined by $R_p(s_0, \dots, s_{p - 1}) = (1 - s_0, \dots, 1 - s_{p - 1})$, which is orientation preserving when $p$ is even and orientation reversing when $p$ is odd.

\begin{myprop}
  The restriction of $F_{i_0, \dots, i_p} \circ R_p$ to $(0, 1)^p$ is an orientation preserving diffeomorphism onto its image.
\end{myprop}
\begin{proof}
  Let $G_{i_0, \dots, i_p} = F_{i_0, \dots, i_p} \circ R_p$.  In coordinates,
  \begin{align*}
    G_{i_0, \dots, i_p}(s_0, \dots, s_{p - 1})) = (1 - s_0)e_{i_0} + \dots + s_0 \dots s_{p - 2}(1 - s_{p - 1}) e_{i_{p - 1}} + s_0\dots s_{p - 1}e_{i_p}
  \end{align*}
  
  We can think of $G_{i_0, \dots, i_p}$ as a map $G_p: [0, 1]^p \to |\Delta^p| \subset \RR^{p + 1}$.  Due to the inverse function theorem, it suffices to show that the differential of $G_p$ is injective on $(0, 1)^p$.  But
  \begin{align*}
    \frac{G^j_p}{\partial s_k} = 0 \quad \text{if $k > j$, and} \quad \frac{G^k_p}{\partial s_k} = - \prod_{l = 0}^{k - 1} s_l
  \end{align*}
  showing that differential has rank $p$ on $(0, 1)^p$ and is therefore injective.   We claim that the determinant of the matrix $(1_{p \times 1} \,  DG_p)$ is positive on $(0, 1)^p$ which implies that $G_{i_0, \dots, i_p}$ is orientation preserving.  The top row of this matrix only has only two nonzero entries.  Therefore, it is sufficient to prove that the two corresponding terms in the cofactor expansion are positive on $(0, 1)^p$, which can be shown by induction on $p$.
\end{proof}

We can now calculate that $\int_{[i_0, \dots, i_p]} \omega_{i_0, \dots, i_p} = 1/p!$ using the identification:
\begin{mylem}\label{lem:integral}
  For any $p$-form $\omega$ on $|\Delta^n|$
  \begin{align*}
    \int_{[i_0, \dots, i_p]}\omega = (-1)^p\varepsilon^{i_p}h^{i_{p - 1}}\dots h^{i_0}(\omega)
  \end{align*}
\end{mylem}
\begin{proof}
Because
\begin{align*}
  \varepsilon^{i_p}h^{i_{p - 1}}\dots h^{i_0}\omega &= \varepsilon^{i_p}\int_{[0,1]^p}(\id_{[0, 1]^{p - 1}} \times \phi^{i_{p - 1}})^{*} \circ \dots \circ (\id_{[0, 1]} \times \phi^{i_1}) \circ  (\phi^{i_0})^{*}\omega\\
                                                    &= \int_{[0, 1]^p}  F_{i_0, \dots, i_p}^{*}\omega  = (-1)^p\int_{[0, 1]^p}  G_{i_0, \dots, i_p}^{*}\omega
\end{align*}
the result follows from the fact that $G_{i_0, \dots, i_p}$ is a parametrization of $[i_0, \dots, i_p]$.
\end{proof}

The proof of this lemma in Getzler's paper \cite{2004math......4003G}, is by arguing by induction on $p$ and assuming that $\omega$ is exact for $p > 0$.  The proof makes use of Stokes' theorem for the induction step.

\subsection{Proof of Main Theorems}
\label{sec:proofs}
The homotopy formula $dh^i + h^id = \varepsilon^i - 1$ is a consequence of the relation $d_{\Delta^k}\pi_{*} = -\pi_{*}d_{\Delta^k}$ and the fundamental theorem of calculus:
\begin{align*}
  \pi_{*}d_{[0, 1]}(\phi^i)^{*}\omega &= (\phi_1^i)^{*}\omega - (\phi_0^i)^{*}\omega\\
  &= \varepsilon^i\omega - \omega.
\end{align*}

\begin{mythm}\label{thm:proof-main-theorem}
  The Dupont map $s$ is a homotopy between $1$ and $WR$.
\end{mythm}
\begin{proof}
  Following Getzler \cite{2004math......4003G}, we compute
\begin{align*}
  [d, s] &= -\sum_{k = 0}^{n - 1}\sum_{i_0 < \dots < i_k}d(\overline{\omega}_{i_0 \dots i_k})h^{i_k}\dots h^{i_0}\\
  &- \sum_{k = 0}^n\sum_{j = 0}^k(-1)^j\sum_{i_0 < \dots < i_k}\overline{\omega}_{i_0 \dots i_k}h^{i_k} \dots [d, h^{i_j}] \dots h^{i_0}
\end{align*}
Using the formula $dh^i + h^id = \varepsilon^i - 1$, the second term above becomes 
\begin{align*}
  1 + \sum_{k = 1}^n\sum_{j = 0}^k(-1)^j\sum_{i_0 < \dots < i_k}\overline{\omega}_{i_0 \dots i_k}h^{i_k} \dots \widehat{h^{i_j}} \dots h^{i_0} - WR.
\end{align*}
The middle term of these three terms is equal to
\begin{align*}
  \sum_{k = 0}^{n - 1}\sum_{i_0 < \dots < i_k}\sum_{i \notin \{i_0, \dots, i_k\}}\overline{\omega}_{ii_0 \dots i_k}h^{i_k}\dots h^{i_0} = \sum_{k = 0}^{n - 1}\sum_{i_0 < \dots < i_k}d(\overline{\omega}_{i_0 \dots i_k})h^{i_k}\dots h^{i_0}
\end{align*}
We conclude that $[d, s] = 1 - WR$.
\end{proof}

To show that $s^2 = 0$, we shall need the identity
\begin{mylem}\label{lem:hiomega}
  If $i \notin \lbrace i_0, \dots, i_k \rbrace$, then
  \begin{align*}
    h^i\overline{\omega}_{i_0, \dots, i_k} = (-1)^k\overline{\omega}_{i_0, \dots, i_k}h^i -  h^i\overline{\omega}_{i_0, \dots, i_k, i}h^i 
  \end{align*}
    Here we are distinguishing between $h^i(\omega)$, which is $h^i$ applied to $\omega$, and $h^i\omega$, which is $h^i$ composed with multiplication by $\omega$.  
\end{mylem}
\begin{proof}
  Getzler \cite{2004math......4003G} observes that we have
  \begin{align*}
    h^i(\omega_{i_0, \dots, i_k}\omega) &= \int_0^1 \frac{ds}{1 - s}(\phi_s^i)^{*}\imath_{E_i}(\omega_{i_0, \dots, i_k} \omega)\\
                                        &= (-1)^k\omega_{i_0, \dots, i_k} \int_0^1ds \, (1 - s)^k (\phi_{s}^i)^{*}\imath_{E_i}\omega 
  \end{align*}
  And the other hand,
  \begin{align*}
    h^i(\omega_{i_0, \dots, i_k, i}h^i(\omega)) &= \int_0^1 \int_0^1 \frac{ds ds'}{(1 - s)(1 - s')} (\phi^i_s)^{*} \imath_{E_i} \omega_{i_0, \dots, i_k, i}(\phi^i_{s'})^{*} \imath_{E_i}\omega\\
    &= (-1)^k\omega_{i_0, \dots, i_k}\int_0^1 \int_0^1 \frac{ds ds' (1 - s)^k}{(1 - s')} (\phi^i_s)^{*} (\phi_{s'})^{*} \imath_{E_i}\omega
  \end{align*}
  Note that $\phi_{1 - s}^i \circ \phi_{1 - s'}^i = \phi_{1 - ss'}^i$.  Upon making the change of variables from $(s, s')$ to $(w, s')$ where $w = ss'$, we see that
  \begin{align*}
    \int_0^1 \int_0^1 \frac{ds ds' (1 - s)^k}{(1 - s')} (\phi^i_s)^{*} (\phi^i_{s'})^{*} &= \int_0^1 \int_0^1 \frac{ds ds' s^k}{s'} (\phi^i_{1 - s})^{*} (\phi^i_{1 - s'})^{*}\\
                                                                                         &= \frac{1}{k + 1} \int_0^1 dw (w^{-1} - w^k)(\phi_{1 - w}^i)^{*}\\
                                                                                         &= \frac{1}{k + 1} \int_0^1\frac{ds}{1 - s} (\phi_s^i)^{*} - \frac{1}{k + 1} \int_0^1ds (1 - s)^k (\phi_s^i)^{*}
  \end{align*}
  establishing that
  \begin{align*}
    (k + 1) h^i\omega_{i_0, \dots, i_k, i}h^i = (-1)^k\omega_{i_0, \dots, i_k}h^i - h^i\omega_{i_0, \dots, i_k}
  \end{align*}
  which establishes the lemma.
\end{proof}

\begin{mythm}\label{thm:dupont-special-main-theorem}
  The \emph{Dupont homotopy} $s$ is a special deformation retraction.
\end{mythm}
\begin{proof}
  Using Lemma \ref{lem:hiomega},
  \begin{align*}
    -h^{i_k}\dots h^{i_0}s &= \sum_{l = 0}^{n - 1}\sum_{\substack{j_0 < \dots < j_l\\ i_0 \in \lbrace j_0, \dots, j_l \rbrace}} h^{i_k}\dots h^{i_1}h^{i_0} \overline{\omega}_{j_0, \dots, j_l} h^{j_l}\dots h^{j_0} \\
    &+ \sum_{l = 0}^{n - 1} \sum_{\substack{j_0 < \dots < j_l\\ i_0 \notin \lbrace j_0, \dots, j_l \rbrace}} h^{i_k}\dots h^{i_1}((-1)^l\overline{\omega}_{j_0, \dots, j_l} - h^{i_0}\overline{\omega}_{j_0, \dots, j_l, i_0})h^{i_0} h^{j_l}\dots h^{j_0} \\
    &= \sum_{l = 0}^{n - 1}(-1)^{(k + 1)l} \sum_{\substack{j_0 < \dots < j_l\\ \lbrace i_0, \dots, i_k \rbrace \cap \lbrace j_0, \dots, j_l \rbrace = \emptyset}} \overline{\omega}_{j_0, \dots, j_l}h^{i_k}\dots h^{i_0}h^{j_l}\dots h^{j_0},
  \end{align*}
  immediately implying that
  \begin{align*}
    s^2 = \sum_{k, l = 0}^{n - 1}(-1)^{(k + 1)l} \sum_{\substack{i_0 < \dots < i_k; j_0 < \dots < j_l\\ \lbrace i_0, \dots, i_k \rbrace \cap \lbrace j_0, \dots, j_l \rbrace = \emptyset}} \overline{\omega}_{i_0, \dots, i_k}\overline{\omega}_{j_0, \dots, j_l}h^{i_k}\dots h^{i_0}h^{j_l}\dots h^{j_0}.
  \end{align*}
  But $h^ih^j + h^jh^i = 0$ for $i \neq j$, so we can interchange the $i$ and $j$ indices in a term to and pick up a factor of $(-1)^{(l+ 1)(k + 1) + lk}$.  The overall sign for an interchanged term is $(-1)^{(k + 1)l}(-1)^{(l+ 1)(k + 1) + lk} = - (-1)^{(l + 1)k}$, so we get pairwise cancellation, implying that $s^2 = 0$.

  It remains to show firstly that $Rs = 0$.  This is a consequence of the fact that the only term in $h^{i_{k - 1}}\dots h^{i_0}s$ that could possibly remain nonzero in the formula above after applying $\varepsilon^{i_k}$ is $\omega_{i_k}h^{i_{k - 1}}\dots h^{i_0}h^{i_k}$, but the identity $\varepsilon^jh^j = 0$ implies this term is also zero.

  And lastly, $sW = 0$ because
  \begin{align*}
    s(\overline{\omega}_{j_0, \dots, j_l}) &= -\sum_{k = 0}^{n - 1}\sum_{i_0 < \dots < i_k} \overline{\omega}_{i_0, \dots, i_k} h^{i_k}\dots h^{i_0}(\overline{\omega}_{j_0, \dots, j_l})
  \end{align*}
  and for each $k$,
  \begin{align*}
    \sum_{i_0 < \dots < i_k}\overline{\omega}_{i_0, \dots, i_k} h^{i_k}\dots h^{i_0} (\overline{\omega}_{j_0, \dots, j_l}) = 0.
  \end{align*}
  This is because $\overline{\omega}_{i_0, \dots, i_k} h^{i_k}\dots h^{i_0}(\overline{\omega}_{j_0, \dots, j_l})$ is equal to
  \begin{align*}
    & k! l! \, i_{\Delta^n}^{*}(\imath_{E}dt_{i_0} \dots dt_{i_k}) i_{\Delta^n}^{*}(\imath_E\imath_{\partial_{i_k}}\dots \imath_{\partial_{i_0}}dt_{j_0}\dots dt_{j_l})\\
    = \, & k! l! \, i_{\Delta^n}^{*}\sum^{l}_{i, j = 0}t_it_j(\imath_{\partial_{i}}dt_{i_0} \dots dt_{i_k}) (\imath_{\partial_j}\imath_{\partial_{i_k}}\dots \imath_{\partial_{i_0}}dt_{j_0}\dots dt_{j_l})
  \end{align*}
  where we have defined $E = \sum_i t_i \partial_{t_i}$.
  The terms in the above sum are zero unless $\lbrace i_0, \dots, i_k \rbrace \subset \lbrace j_0, \dots, j_l \rbrace$, $i \in \lbrace i_0, \dots, i_k \rbrace$, and $j \in  \lbrace j_0, \dots, j_l \rbrace \setminus \lbrace i_0, \dots, i_k\rbrace$.  So in the sum over $i_0 < \dots < i_k$, we can cancel the term indexed by $i_0 < \dots < i_k$ with its pair: the term indexed by $i'_0 < \dots < i'_k$ which is given by removing $i$ adding $j$.
\end{proof}

\subsection{Globalizing the Construction}
\label{sec:global_const}
It is important to clarify what is meant by a piecewise smooth differential form.  To specify a piecewise smooth form, we specify a smooth form $\omega_T$ for each simplex $T$ with the compatibility condition that the pullback of $\omega_T$ to a face $T'$ of $T$ is $\omega_{T'}$.

The general statement for a triangulated manifold follows as a consequence of Theorems \ref{thm:proof-main-theorem} \& \ref{thm:dupont-special-main-theorem} as we now show:

\begin{mycor}
  The \emph{Dupont homotopy} $s$ is a well-defined (special) deformation retraction of the piecewise smooth differential forms $\Omega^{\bullet}_{\Xi}(M)$ on a triangulated smooth manifold $M \cong |\Xi|$ onto the simplicial cochains $C^{\bullet}(\Xi)$.  
\begin{align*}
  \myDR{(C^{\bullet}(\Xi), d)}{(\Omega^{\bullet}_{\Xi}(M), d)}{W}{R}{s}
\end{align*}
where $W$ is the Whitney map and $R$ is the integration map, which are also well-defined.
\end{mycor}
\begin{proof}
  Firstly, we need on $\Delta^n$ that $W$, $R$ and $s$ are equivariant under the action of the symmetric group.  This means that the maps do not depend on the ordering that we choose for the vertices of $\Delta^n$.  This is necessary because a triangulation of a manifold does not come with an ordering of the vertices of its simplices.  For $W$ and $R$ this follows directly from the definition.  Turning to $s$ for a general permutation $\sigma \in S_n$ let $\tau_{\sigma}: |\Delta^n| \to |\Delta^n|$ be the induced map $\tau_{\sigma}(t_0, \dots, t_n) = (t_{\sigma(0)}, \dots, t_{\sigma(n)})$.  We have
  \begin{align*}
    \tau_{\sigma}^{*} s &=  -\sum_{k = 0}^{n - 1} \sum_{i_0 < \dots < i_k} \tau_{\sigma}^{*}\overline{\omega}_{i_0, \dots, i_k} \tau_{\sigma}^{*} h^{i_k}\dots h^{i_0} \\
                 &= -\sum_{k = 0}^{n - 1} \sum_{i_0 < \dots < i_k} \overline{\omega}_{\sigma(i_0), \dots, \sigma(i_k)}  h^{\sigma(i_k)}\dots h^{\sigma(i_0)} \tau_{\sigma}^{*}\\
    &= s \tau_{\sigma}^{*}.                   
  \end{align*}
Here we have used the fact that $\tau_{\sigma}^{*}h^i = h^{\sigma(i)}\tau_{\sigma}^{*}$ which follows from the identity $\tau_{\sigma} \circ \phi^{\sigma(i)} = \phi^i \circ (\id_{[0, 1]} \times \tau_{\sigma})$.
  
  It now suffices to show that $W$, $R$ and $s$ commute with pullback by the face maps $\epsilon_i: \Delta^{n - 1} \to \Delta^n$ and their geometric realizations $|\epsilon_i|: |\Delta^{n - 1}| \to |\Delta^n|$ for $i = 0, \dots, n$.  For $W$ and $R$ this follows directly from the definition.  For $s$ we use that
  \begin{align*}
    |\epsilon_i| \circ \phi^j =
    \begin{cases}
      \phi^j \circ (\id_{[0, 1]} \times |\epsilon_i|) & \text{if $i > j$}\\ 
      \phi^{j + 1} \circ (\id_{[0, 1]} \times |\epsilon_i|) & \text{if $i \leq j$}
    \end{cases}
  \end{align*}
  This implies that
  \begin{align*}
    h^j|\epsilon_i|^{*} =
    \begin{cases}
      |\epsilon_i|^{*}h^j & \text{if $i > j$}\\
      |\epsilon_i|^{*}h^{j + 1} & \text{if $i \leq j$}
    \end{cases}
  \end{align*}
  and therefore
  \begin{align*}
    |\epsilon_i|^{*}s &=  -\sum_{k = 0}^{n - 1} \sum_{\substack{i_0 < \dots < i_k \\ i \notin \{i_0, \dots, i_k \}}} |\epsilon_i|^{*}(\overline{\omega}_{i_0, \dots, i_k})|\epsilon_i|^{*} h^{i_k}\dots h^{i_0} \\
    &= -\sum_{k = 0}^{n - 1} \sum_{\substack{i_0 < \dots < i_j < i \\ i < i_{j + 1} < \dots < i_k}} \overline{\omega}_{i_0, \dots, i_j, i_{j + 1} - 1, \dots, i_k - 1} h^{i_k - 1}\dots h^{i_{j + 1} - 1} h^{i_j}\dots h^{i_0} |\epsilon_i|^{*} \\
    &= s|\epsilon_i|^{*}
  \end{align*}
\end{proof}

\section{Cubical Dupont Homotopy Formula}
\label{sec:cub_dup_hom}
The Dupont homotopy formula for cubical forms can be constructed from the Dupont homotopy formula on the $1$-simplex $[0, 1]$ through the tensor product construction.

Let $t$ be the natural coordinate on $|\Delta^1| \cong [0, 1]$.  The degree $0$ Whitney forms on $|\Delta^1|$ are $\omega_0 = 1 - t$ and $\omega_1 = t$.  The degree $1$ Whitney form is
\begin{align*}
  \omega_{01} = (1 - t)dt - td(1 - t) = dt
\end{align*}
Given a form $\omega = f(t) + g(t)\, dt$ on $|\Delta^1| = [0, 1]$, 
\begin{align*}
  RW(\omega) &= \omega_0 R_0\omega + \omega_1 R_1\omega + \omega_{01} R_{0, 1}\omega\\
  &= (1 - t) f(0) + tf(1) + dt \int_0^1g(t)\, dt \
\end{align*}
and the Dupont homotopy is given by
\begin{align*}
  s(\omega) &= -\omega_0 h^0(\omega) - \omega_1 h^1(\omega)\\
            &= (1 - t) t\int_0^1g((1 - s)t) \, ds + t(t - 1) \int_0^1 g((1 - s)t + s)\, ds\\
            &=  (1 - t) \int_0^tg(u)\, du - t \int_t^1 g(u) \, du \\
            &= \int_0^tg(u)\, du - t \int_0^1 g(u) \, du              
\end{align*}

Let us recall the tensor construction.  Note that $C^{\bullet}(\square^n) = C^{\bullet}(\Delta^1) \otimes \dots \otimes C^{\bullet}(\Delta^1)$ and $\Omega^{\bullet}(|\square^n|) = \Omega^{\bullet}(|\Delta^1|) \widehat{\otimes} \dots \widehat{\otimes} \Omega^{\bullet}(|\Delta^1|)$ where $\widehat{\otimes}$ denotes the completed projective tensor product.  Due to the continuity of $R: \Omega^{\bullet}(|\Delta^1|) \to C^{\bullet}(\Delta^1)$, $WR: \Omega^{\bullet}(|\Delta^1|) \to \Omega^{\bullet}(|\Delta^1|)$ and $s: \Omega^{\bullet}(|\Delta^1|) \to \Omega^{\bullet}(|\Delta^1|)$, the following definitions make sense:
\begin{mydef}
  We define the \emph{integration map} $\mathbf{R}: \Omega^{\bullet}(|\square^n|) \to C^{\bullet}(\square^n)$ by $\mathbf{R} = R \otimes \dots \otimes R$ and the \emph{cubical Whitney map} $\mathbf{W}: C^{\bullet}(\square^n) \to \Omega^{\bullet}(|\square^n|)$ by $\mathbf{W} = W \otimes \dots \otimes W$.  Define
  \begin{align*}
    \mathbf{s}_0 =  \sum_{j = 1}^n \underbrace{1 \otimes \dots \otimes 1}_{j - 1} \otimes s \otimes WR \otimes \dots \otimes WR
  \end{align*}
and define the \emph{cubical Dupont homotopy} $\mathbf{s}$ as the symmetrization of $\mathbf{s}_0$.  That is if $\tau_{\sigma}:\Omega^{\bullet}(|\square^n|) \to \Omega^{\bullet}(|\square^n|)$ is the induced linear map coming from the permutation $\sigma \in S_n$, we have
\begin{align*}
  \mathbf{s} &= \frac{1}{n!} \sum_{\sigma \in S_n} \tau_{\sigma} \circ \mathbf{s}_0 \circ \tau_{\sigma}^{-1}\\
             &= \frac{1}{n!}\sum_{\mathbf{\epsilon}} C_{|\epsilon|, n} \sum_{j = 1}^n (WR)^{\epsilon_1} \otimes \dots \otimes (WR)^{\epsilon_{j - 1}} \otimes s \otimes (WR)^{\epsilon_j}\otimes \dots \otimes (WR)^{\epsilon_{n - 1}}\\
  &= \frac{1}{n!}\sum_{\mathbf{\epsilon}} C_{|\epsilon|, n} \boldsymbol{\psi}_{\epsilon}
\end{align*}
where $C_{|\epsilon|, n} = |\epsilon|!(n - 1 - |\epsilon|)!$ and the outer sum is over $\epsilon$ with $\epsilon_k = 0, 1$.  
\end{mydef}

\begin{mythm}
  The \emph{cubical Dupont homotopy} $\textbf{s}$ is a (special) deformation retraction of the differential forms $\Omega^{\bullet}(|\square^n|)$ on the $n$-cube onto the cubical cochains $C^{\bullet}(\square^n)$.  
\begin{align*}
  \myDR{(C^{\bullet}(\square^n), d)}{(\Omega^{\bullet}(|\square^n|), d)}{\mathbf{W}}{\mathbf{R}}{\mathbf{s}}
\end{align*}
where $\textbf{W}$ is the Whitney map and $\textbf{R}$ is the integration map.
\end{mythm}
\begin{proof}  
  The theorem also holds replacing $\mathbf{s}$ with $\mathbf{s}_0$.  The reason for working with $\mathbf{s}$ rather than $\mathbf{s}_0$ is to be able to pass to cubulated manifolds where there is no fixed identification of the $n$-cube as an ordered product of $1$-simplices.

  Because $d$ commutes with $\tau_{\sigma}$ for any permutation $\sigma \in S_n$, to show that $\mathbf{s}$ is a deformation retraction, it suffices to show that $\mathbf{s}_0$ is a deformation retraction.  But
  \begin{align*}
    d \mathbf{s}_0 +  \mathbf{s}_0 d  &= \sum_{j =  1}^n 1^{\otimes (j - 1)} \otimes (ds + sd) \otimes (WR)^{\otimes(n - j)} \\
                                  &= \mathbf{1} - \mathbf{WR}
  \end{align*}

  It is clear that $\mathbf{s} \mathbf{W} = 0$ and $\mathbf{R} \mathbf{s} = 0$.  Lastly $\mathbf{s}^2 = 0$ follows from the relation $\boldsymbol{\psi}_{\epsilon} \boldsymbol{\psi}_{\epsilon'} = -\boldsymbol{\psi}_{\epsilon'} \boldsymbol{\psi}_{\epsilon}$.
\end{proof}

\begin{mycor}
  The \emph{cubical Dupont homotopy formula} gives a (special) deformation retraction of the piecewise smooth differential forms $\Omega_{\Xi}^{\bullet}(M)$ on a cubulated manifold $M$ onto the cubical cochains $C^{\bullet}(\Xi)$.  
\begin{align*}
  \myDR{(C^{\bullet}(\Xi), d)}{(\Omega_{\Xi}^{\bullet}(M), d)}{W}{R}{s}
\end{align*}
where $W$ is the cubical Whitney map, $R$ is the cubical integration map, and $s$ is the cubical Dupont homotopy. 
\end{mycor}

\section{Stellar Subdivision Formulas}
\label{sec:stellar_subdiv_form}
We begin with the statement in one dimension for simplicity.  Let $\star \Delta^1$ denote the stellar subdivision of the $1$-simplex $\Delta^1$.  That is $\star \Delta^1$ is the simplicial complex with vertices $e_{\star}, e_0, e_1$ and edges $[e_0, e_{\star}]$ and $[e_{\star}, e_1]$.  We visualize $e_{\star}$ as lying at the barycenter of the $1$-simplex $[e_0, e_1]$ of $\Delta^1$.
\begin{center}
  \begin{tikzpicture}
    \tikzstyle{vertex} = [circle, fill=black, inner sep=1pt]
    \node[vertex, label=below:$e_0$] (A1) at (0, 1) {};
    \node[vertex, label=below:$e_1$] (B1) at (2, 1) {};
    \draw (A1) -- (B1);

    \node[vertex, label=below:$e_0$] (A) at (4, 1) {};
    \node[vertex, label=below:$e_\star$] (B) at (5, 1) {};
    \node[vertex, label=below:$e_1$] (C) at (6, 1) {};
    \draw (A) -- (B) -- (C);
  \end{tikzpicture}
\end{center}

There is a natural inclusion map $\imath_{\star}: C_{\bullet}(\Delta^1) \to C_{\bullet}(\star \Delta^1)$ defined by $e_0 \mapsto e_0$, $e_1 \mapsto e_1$ and $[e_0, e_1] \mapsto [e_0, e_{\star}] + [e_{\star}, e_1]$, which is a chain map.  There is a natural projection map $p_{\star}: C_{\bullet}(\star \Delta^1) \to C_{\bullet}(\Delta^1)$ defined by $e_{\star} \mapsto \frac{1}{2}(e_0 + e_1)$, $e_0 \mapsto e_0$, $e_1 \mapsto e_1$ and $[e_0, e_{\star}] \mapsto \frac{1}{2}[e_0, e_1]$ and $[e_{\star}, e_1] \mapsto \frac{1}{2}[e_0, e_1]$, which is also a chain map.  We have $p_{\star}i_{\star} = 1$ and we would like to find a homotopy $a_{\star}$ between the identity $1$ and $i_{\star}p_{\star}$.  We define $a_{\star}$ by $e_{\star} \mapsto \frac{1}{2}([e_0, e_{\star}] - [e_{\star}, e_1])$, $e_0 \mapsto 0$, and $e_1 \mapsto 0$.  Then $\partial a_{\star} + a_{\star}\partial = 1 - i_{\star}p_{\star}$.  Furthermore, we have $a_{\star}i_{\star} = 0$, $p_{\star}a_{\star} = 0$, and $a_{\star}^2 = 0$.

The dual deformation retraction in one dimension has inclusion $\imath^{\star}: C^{\bullet}(\Delta^1) \to C^{\bullet}(\star \Delta^1)$ defined by $\widehat{e}_0 \mapsto \widehat{e}_0 + \frac{1}{2} \widehat{e}_{\star}$, $\widehat{e}_1 \mapsto \widehat{e}_1 + \frac{1}{2} \widehat{e}_{\star}$, $\widehat{[e_0, e_1]} \mapsto \frac{1}{2}( \widehat{[e_0, e_{\star}]} + \widehat{[e_{\star}, e_1]})$.  It has projection $p^{\star}: C^{\bullet}(\star\Delta^1) \to C^{\bullet}(\Delta^1)$ defined by $\widehat{e}_0 \mapsto \widehat{e}_0$, $\widehat{e}_1 \mapsto \widehat{e}_1$, $\widehat{e}_{\star} \mapsto 0$, $\widehat{[e_0, e_{\star}]} \mapsto \widehat{[e_0, e_1]}$ and $\widehat{[e_{\star}, e_1]} \mapsto \widehat{[e_0, e_1]}$.  Lastly, the homotopy $a^{\star}$ defined by $\widehat{[e_0, e_{\star}]} \mapsto \frac{1}{2}e_{\star}$ and $\widehat{[e_{\star}, e_1]} \mapsto -\frac{1}{2}e_{\star}$.

Generalizing now to the $n$-simplex:
\begin{mydef}\cite{Kozlov_2008}
  For $k \leq n$ and $0 \leq i_0 < \dots < i_k \leq n$, we define the \emph{stellar subdivison} $\star_{i_0, \dots, i_k} \Delta^n$ for $k \leq n$.  This is a simplicial complex having vertex $e_{\star}$ as well as vertices $e_0, \dots, e_n$.  For $J \not \supset I$, where $J = \{j_0, \dots, j_l\}$ and $I = \{i_0, \dots, i_k\}$, we include all simplices of the form $[e_{j_0}, \dots, e_{j_l}]$ and $[e_{\star}, e_{j_0}, \dots, e_{j_l}]$.  When $k = n$ we shall sometimes just write $\star \Delta^n$.  We visualize $e_{\star}$ as lying at the barycenter of simplex $[e_{i_0}, \dots, e_{i_k}]$.
\end{mydef}

\begin{center}
  \begin{tikzpicture}
    \tikzstyle{vertex} = [circle, fill=black, inner sep=1pt]
    \node[vertex, label=below:$e_0$] (A1) at (0, 1) {};
    \node[vertex, label=below:$e_2$] (B1) at (2, 1) {};
    \node[vertex, label=right:$e_1$] (C1) at (1, 2.5) {};
    \draw (A1) -- (B1) -- (C1) -- (A1);

    \node[vertex, label=below:$e_0$] (A1) at (4, 1) {};
    \node[vertex, label=below:$e_2$] (B1) at (6, 1) {};
    \node[vertex, label=right:$e_1$] (C1) at (5, 2.5) {};
    \node[vertex, label=below:$e_\star$] (D1) at (5, 1.5) {};
    \draw (A1) -- (B1) -- (C1) -- (A1);
    \draw (A1) -- (D1);
    \draw (B1) -- (D1);
    \draw (C1) -- (D1);
  \end{tikzpicture}
\end{center}

\begin{mydef}
  We define the \emph{stellar welding inclusion map} $i_{\star}: C_{\bullet}(\Delta^n) \to C_{\bullet}(\star_{i_0, \dots, i_k}\Delta^n)$ by 
  \begin{align*}
    i_{\star}[e_{j_0}, \dots, e_{j_l}] =
    \begin{cases}
      [e_{j_0}, \dots, e_{j_l}] & \text{ for $J \not \supset I$} \\
      \sum_{j_i \in I} (-1)^i[e_{\star}, e_{j_0}, \dots, \widehat{e_{j_i}}, \dots,  e_{j_l}] & \text{ for $J \supset I$} 
    \end{cases}
  \end{align*}
\end{mydef}

We quickly verify that $i_{\star}$ is a chain map since for $J \supset I$
\begin{align*}
  \partial i_{\star} [e_{j_0}, \dots, e_{j_l}] &= i_{\star}\sum_{j_i \in I} (-1)^i[e_{j_0}, \dots, \widehat{e_{j_i}}, \dots,  e_{j_l}] \\
                                               &+ i_{\star} \sum_{j_i \not \in I} (-1)^i[e_{j_0}, \dots, \widehat{e_{j_i}}, \dots,  e_{j_l}] \\
  &= i_{\star} \partial  [e_{j_0}, \dots, e_{j_l}]
\end{align*}

\begin{mydef}
  Define the \emph{stellar welding projection map} $p_{\star}: C_{\bullet}(\star_{i_0, \dots, i_k} \Delta^n) \to C_{\bullet}(\Delta^n)$ for $J \not \supset I$, by
  \begin{align*}
    p_{\star} [e_{j_0}, \dots, e_{j_l}] = [e_{j_0}, \dots, e_{j_l}]
  \end{align*}
  and
  \begin{align*}
    p_{\star} [e_{\star}, e_{j_0}, \dots, e_{j_l}] &= \frac{1}{k + 1}\sum_{\alpha \in I \setminus J} [e_\alpha, e_{j_0}, \dots, e_{j_l}]
  \end{align*}

\end{mydef}
This is a chain map because
\begin{align*}
  \partial p_{\star} [e_{\star}, e_{j_0}, \dots, e_{j_l}] &=  \frac{|I \setminus J|}{k + 1}[e_{j_0}, \dots, e_{j_l}] - \frac{1}{k + 1} \sum_{\alpha \in I \setminus J} \sum_{i = 0}^l \, (-1)^{j_i} [e_\alpha, e_{j_0}, \dots, \widehat{e_{j_i}}, \dots, e_{j_l}]\\
                                                          &= [e_{j_0}, \dots, e_{j_l}] - \frac{1}{k + 1} \sum_{i = 0}^l \sum_{\alpha \in I \setminus (J \setminus \{j_i\})}  (-1)^{j_i}[e_\alpha, e_{j_0}, \dots, \widehat{e_{j_i}}, \dots, e_{j_l}]\\ 
                                                          &= p_{\star} \partial [e_{\star}, e_{j_0}, \dots, e_{j_l}] 
\end{align*}

\begin{mydef}
  Define the \emph{stellar welding homotopy} by $a_{\star}: C_{\bullet}(\star_{i_0, \dots, i_k} \Delta^n) \to C_{\bullet + 1}(\star_{i_0, \dots, i_k} \Delta^n)$ for $J \not \supset I$ by
  \begin{align*}
    a_{\star}[e_{j_0}, \dots, e_{j_l}] = 0 
  \end{align*}
  \begin{align*}
    a_{\star}[e_{\star}, e_{j_0}, \dots, e_{j_l}] =
    \begin{cases}
      0 & \text{ if $|I \setminus J| = 1$}\\
      - \frac{1}{k + 1}\sum_{\alpha \in I \setminus J}[e_{\star}, e_{\alpha}, e_{j_0}, \dots, e_{j_l}] & \text{ otherwise}
    \end{cases}
  \end{align*}
\end{mydef}

\begin{mythm}\label{sec:stell-subd-chains}
  \begin{align*}
    \partial a_{\star} + a_{\star} \partial = 1 - i_{\star}p_{\star}
  \end{align*}
  and $a_{\star}^2 = 0$, $a_{\star}i_{\star} = 0$, and $p_{\star}a_{\star} = 0$.
\end{mythm}
\begin{proof}
  We compute for $J \not \supset I$ 
\begin{align*}
  i_{\star} p_{\star}[e_{j_0}, \dots, e_{j_l}] = [e_{j_0}, \dots, e_{j_l}],
\end{align*}
and
\begin{align*}
  (\partial a_{\star} + a_{\star}\partial)([e_{j_0}, \dots, e_{j_l}]) = 0
\end{align*}

For $\sigma = [e_{\star}, e_{j_0}, \dots, e_{j_l}]$ with $I \setminus J = \{i_m\}$,
\begin{align*}
  i_{\star} p_{\star}\sigma = \frac{1}{k + 1}[e_{\star}, e_{j_0}, \dots, e_{j_l}] - \frac{1}{k + 1}\sum_{j_i \in I} (-1)^i[e_{\star}, e_{i_m}, e_{j_0}, \dots, \widehat{e_{j_i}}, \dots, e_{j_l}] 
\end{align*}
and 
\begin{align*}
  (\partial a_{\star} + a_{\star}\partial)([e_{\star}, e_{j_0}, \dots, e_{j_l}]) &= -\sum_{i = 0}^l (-1)^i a_{\star} [e_{\star}, e_{j_0}, \dots, \widehat{e_{j_i}}, \dots, e_{j_l}] \\
                                                                                 &= \frac{k}{k + 1} [e_{\star}, e_{j_0}, \dots, e_{j_l}] \\
                                                                                 &+ \frac{1}{k + 1}\sum_{j_i \in I} (-1)^i[e_{\star}, e_{i_m}, e_{j_0}, \dots, \widehat{e_{j_i}}, \dots, e_{j_l}]\\
  &= (1 - i_{\star}p_{\star})([e_{\star}, e_{j_0}, \dots, e_{j_l}])
\end{align*}

For $\sigma = [e_{\star}, e_{j_0}, \dots, e_{j_l}]$ with $|I \setminus J| \geq 2$, we have
\begin{align*}
  i_{\star} p_{\star}\sigma =  \frac{1}{k + 1}\sum_{\alpha \in I \setminus J} [e_\alpha, e_{j_0}, \dots, e_{j_l}]
\end{align*}
and 
\begin{align*}
  (\partial a_{\star} + a_{\star}\partial)([e_{\star}, e_{j_0}, \dots, e_{j_l}]) &= -\sum_{i = 0}^l (-1)^i a_{\star} [e_{\star}, e_{j_0}, \dots, \widehat{e_{j_i}}, \dots, e_{j_l}] \\
            &- \frac{1}{k + 1}\sum_{\alpha \in I \setminus J}\partial [e_{\star}, e_{\alpha}, e_{j_0}, \dots, e_{j_l}] \\
  &= (1 - i_{\star}p_{\star})([e_{\star}, e_{j_0}, \dots, e_{j_l}])
\end{align*}

It is clear from the definitions that $a_{\star}i_{\star} = 0$.  For $|I \setminus J| \geq 2$, we have
\begin{align*}
  p_{\star}a_{\star}[e_{\star}, e_{j_0}, \dots, e_{j_l}] = -\frac{1}{(k + 1)^2}\sum_{\alpha \in I \setminus J} \sum_{\alpha' \in I \setminus (J \cup \{\alpha\})}  [e_{\alpha'}, e_{\alpha}, e_{j_0}, \dots, e_{j_l}] = 0
\end{align*}
and if $|I \setminus J| \geq 3$,
\begin{align*}
  a_{\star}^2[e_{\star}, e_{j_0}, \dots, e_{j_l}] = \frac{1}{(k + 1)^2}\sum_{\alpha \in I \setminus J} \sum_{\alpha' \in I \setminus (J \cup \{\alpha\})} [e_{\star}, e_{\alpha'}, e_{\alpha}, e_{j_0}, \dots, e_{j_l}] = 0
\end{align*}
\end{proof}

\begin{mydef}
  We define the \emph{stellar welding inclusion map} $i^{\star}: C^{\bullet}(\Delta^n) \to C^{\bullet}(\star_{i_0, \dots, i_k}\Delta^n)$ by
  \begin{align*}
    i^{\star}\reallywidehat{[e_{j_0}, \dots, e_{j_l}]} = \reallywidehat{[e_{j_0}, \dots, e_{j_l}]} + \frac{1}{k + 1} \sum_{j_i \in I} (-1)^i  \reallywidehat{[e_{\star}, e_{j_0}, \dots, \widehat{e_{j_i}}, \dots, e_{j_l}]}
  \end{align*}
  for $J \not \supset I$ and
  \begin{align*}
    i^{\star} \reallywidehat{[e_{j_0}, \dots, e_{j_l}]} = \frac{1}{k + 1} \sum_{j_i \in I} (-1)^i \reallywidehat{[e_{\star}, e_{j_0}, \dots, \widehat{e_{j_i}}, \dots, e_{j_l}]}
  \end{align*}
  for $J \supset I$.
\end{mydef}

\begin{mydef}
  We define the \emph{stellar welding projection map} $p^{\star}: C^{\bullet}(\star_{i_0, \dots, i_k} \Delta^n) \to C^{\bullet}(\Delta^n)$ by 
  \begin{align*}
    p^{\star} \reallywidehat{[e_{j_0}, \dots, e_{j_l}]} = \reallywidehat{[e_{j_0}, \dots, e_{j_l}]}
  \end{align*}
  and 
  \begin{align*}
    p^{\star} \reallywidehat{[e_{\star}, e_{j_0}, \dots, e_{j_l}]} = 
    \begin{cases}
      \reallywidehat{[e_{i_m}, e_{j_0}, \dots, e_{j_l}]} & \text{ if $I \setminus J = \{i_m\}$}\\
      0 & \text{ otherwise}
    \end{cases}
  \end{align*}
  for $J \not \supset I$.
\end{mydef}
\begin{mydef}
  We define the \emph{stellar welding homotopy} $a^{\star}: C^{\bullet}(\star_{i_0, \dots, i_k} \Delta^n) \to C^{\bullet - 1}(\star_{i_0, \dots, i_k} \Delta^n)$ by 
  \begin{align*}
    a^{\star}\reallywidehat{[e_{j_0}, \dots, e_{j_l}]} = 0
  \end{align*}
and
  \begin{align*}
    a^{\star}\reallywidehat{[e_{\star}, e_{j_0}, \dots, e_{j_l}]} = -\frac{1}{k + 1} \sum_{j_i \in I} (-1)^i \reallywidehat{[e_{\star}, e_{j_0}, \dots, \reallywidehat{e_{j_i}}, \dots, e_{j_l}]} 
  \end{align*}
   for $J \not \supset I$
\end{mydef}

By dualizing the result of Theorem \ref{sec:stell-subd-chains}, we have
\begin{mythm}\label{thm:stell-subd-cochains}
  The maps above define a special deformation retraction
  \begin{align*}
    \myDR{(C^{\bullet}(\Delta^n), d)}{(C^{\bullet}(\star_{i_0, \dots, i_k} \Delta^n), d)}{i^{\star}}{p^{\star}}{a^{\star}}
  \end{align*}
  We call it the \emph{stellar welding deformation retraction (on cochains)}. 
\end{mythm}

\section{Cubical Stellar Subdivision Formulas}
\label{sec:cub_stell_subdiv}
We define the cubical subdivision $\star_{i_1, \dots, i_k} \square^n$ for $\{i_1, \dots, i_k\} \subset \{1, \dots, n\}$, to be the product $I_1 \times \dots \times I_n$ where $I_i = \star \Delta^1$ if $i \in \{i_1, \dots, i_k\}$ and $I_i = \Delta^1$ otherwise.  To simplify notation, we shall always assume in what follows that $i_j = j$ so that
\begin{align*}
  \star_{i_1, \dots, i_k} \square^n = \star_{1, \dots, k} \square^n = (\star\Delta^1)^{\times k} \times (\Delta^1)^{\times (n - k)}.
\end{align*}
The formulas for the cubical welding deformation retraction that we will present can be recovered for general $\{i_1, \dots, i_k\}$ by composing with an appropriate permutation.
\begin{mydef}
  Define $\mathbf{i}_{\star} = i_{\star}^{\otimes k} \otimes 1^{\otimes (n - k)}$, $\mathbf{p}_{\star} = p_{\star}^{\otimes k} \otimes 1^{\otimes (n - k)}$ and
  \begin{align*}
    \mathbf{a}_{\star} = \frac{1}{k!} \sum_{\sigma \in S_k} \tau_{\sigma \times 1^{n - k}} \circ \left(\sum_{j = 1}^k 1^{\otimes (j - 1)} \otimes a_{\star} \otimes (i_{\star}p_{\star})^{\otimes (k - j)} \otimes 1^{\otimes (n - k)} \right)\circ \tau_{\sigma \times 1^{n - k}}^{-1}
  \end{align*}
  respectively the \emph{cubical welding inclusion, projection and homotopy (on cubical chains)} .  Define $\mathbf{i}^{\star} = (i^{\star})^{\otimes k} \otimes 1^{\otimes (n - k)}$, $\mathbf{p}^{\star} = (p^{\star})^{\otimes k} \otimes 1^{\otimes (n - k)}$ and
  \begin{align*}
    \mathbf{a}^{\star} = \frac{1}{k!} \sum_{\sigma \in S_k} \tau_{\sigma \times 1^{n - k}} \circ \left(\sum_{j = 1}^k 1^{\otimes (j - 1)} \otimes a^{\star} \otimes (i^{\star}p^{\star})^{\otimes (k - j)} \otimes 1^{\otimes (n - k)} \right)\circ \tau_{\sigma \times 1^{n - k}}^{-1}
  \end{align*}
  respectively the \emph{cubical stellar subdivision inclusion, projection and homotopy (on cubical cochains)}. 
\end{mydef}

\begin{mythm}
  The above cubical welding maps (on cubical chains) are a special deformation retraction
  whose dual is the special deformation retraction
  \begin{align*}
    \myDR{(C^{\bullet}(\square^n), d)}{(C^{\bullet}(\star_{1, \dots, k} \square^n), d)}{\mathbf{i}^{\star}}{\mathbf{p}^{\star}}{\mathbf{a}^{\star}}
  \end{align*}
  We call this the \emph{cubical welding deformation retraction (on cubical cochains)}.
\end{mythm}

Lastly, it is worth mentioning that cubical subdivision on a cubulated manifold has a different flavor than simplicial stellar subdivision on a triangulated manifold.  Simplicial stellar subdivision is a local operation that exists for any choice of a $k$-simplex in the triangulation.  However, for cubical subdivision, one must specify a collection of $k$-cubes in the cubulation in a compatible way.

\section{Compatibility of the DHF and Stellar Subdivision I}
\label{sec:comp-stell-dup}
We begin by examining the $1$-simplex where the formulas and proofs are much simpler.  There is a deformation retraction
\begin{align*}
  \myDR{(C^{\bullet}(\star\Delta^1), d)}{(\Omega^{\bullet}(\left|\star\Delta^1\right|), d)}{\accentset{\star}{W}}{\accentset{\star}{R}}{\accentset{\star}{s}}
\end{align*}
Let $\chi_{[a, b]}$ be the characteristic function for the interval $[a, b]$.  We identify $\left|\star\Delta^1\right|$ with the interval $[0, 1]$.  Explicitly, the Whitney forms are $\omega_{\star} = 2t\chi_{[0, 1/2]} + (2 - 2t)\chi_{[1/2, 1]}$, $\omega_0 = (1 - 2t)\chi_{[0, 1/2]}$, $\omega_1 = (2t - 1)\chi_{[1/2, 1]} $, $\omega_{0,\star} = \chi_{[0, 1/2]} 2 dt $, and $\omega_{\star, 1} = \chi_{[1/2, 1]} 2 dt$.  The Dupont homotopy is given for $\omega(t) = g(t)\, dt$ by
\begin{align*}
  \myoset{\star}{s}(\omega) &= -\omega_0h^0(\omega) - \omega_{\star}h^{\star}(\omega) - \omega_1h^1(\omega)\\
                                &=\left[\int_0^tg(u) \, du - 2t \int_0^{1/2} g(u)\, du\right] \chi_{[0, 1/2]}\\
                                &+ \left[\int_{1/2}^tg(u) \, du - 2(t - 1/2) \int_{1/2}^1 g(u)\, du\right] \chi_{[1/2, 1]} 
\end{align*}

\begin{mythm}\label{thm:comp-dhf-stell-1}
The deformation retraction of $\Omega^{\bullet}(\left|\star \Delta^1\right|)$ onto $C^{\bullet}(\star \Delta^1)$ can be composed with the stellar welding deformation retraction to get a new deformation retraction
\begin{align*}
  \myDR{(C^{\bullet}(\Delta^1), d)}{(\Omega^{\bullet}(\left|\star\Delta^1\right|), d)}{\accentset{\star}{W} i^{\star}}{p^{\star} \accentset{\star}{R}}{\accentset{\star}{s} + \accentset{\star}{W} a^{\star} \accentset{\star}{R}}
\end{align*}
Comparing this to the deformation retraction
\begin{align*}
  \myDR{(C^{\bullet}(\Delta^1), d)}{(\Omega^{\bullet}(\left|\Delta^1\right|), d)}{W}{R}{s}
\end{align*}
we find that $(p^{\star} \myoset{\star}{R})\mid_{\Omega^{\bullet}(\left|\Delta^1\right|)} = R$, $\myoset{\star}{W} i^{\star} = W$, and $(\myoset{\star}{W} a^{\star} \myoset{\star}{R} + \myoset{\star}{s})\mid_{\Omega^{\bullet}(\left|\Delta^1\right|)} = s$
\end{mythm}
\begin{proof}
We have
\begin{align*}
  p^{\star} \myoset{\star}{R}(f(t)) &= p^{\star}(f(0) \widehat{e_0} + f(1/2) \widehat{e_\star} + f(1) \widehat{e_1}) \\
                                              &= f(0) \widehat{e_0} + f(1) \widehat{e_1}\\
                                              &= R(f(t))
\end{align*}
and
\begin{align*}
  p^{\star} \myoset{\star}{R}(g(t)\, dt) &= p^{\star}\left(\widehat{[e_0, e_\star]}\int_0^{1/2}g(t)\, dt + \widehat{[e_\star, e_1]}\int_{1/2}^1g(t)\, dt\right)\\
                                                   &= \widehat{[e_0, e_1]}\int_0^1g(t)\, dt\\
                                                   &= R( g(t) \, dt)
\end{align*}
That is $p^{\star} \myoset{\star}{R} = R$.

Secondly,
\begin{align*}
  \myoset{\star}{W} i^{\star}(\widehat{e_0}) &= \omega_0 + \frac{1}{2}\omega_\star =  1 - t = W(\widehat{e_0})
\end{align*}
and
\begin{align*}
  \myoset{\star}{W} i^{\star}(\widehat{e_1}) &= \frac{1}{2}\omega_{\star} + \omega_2 = t = W(\widehat{e_1})
\end{align*}
and lastly
\begin{align*}
\myoset{\star}{W} i^{\star}(\widehat{[e_0, e_1]}) = \frac{1}{2}W(\widehat{[e_0, e_\star]} + \widehat{[e_\star, e_1]}) = dt.
\end{align*}
That is, $\myoset{\star}{W} i^{\star} = W$.

Finally, we have that
\begin{align*}
  \myoset{\star}{W} a^{\star} \myoset{\star}{R}(g(t) \, dt) &= \myoset{\star}{W}\left(\frac{\widehat{e_\star}}{2}\int_0^{1/2} g(t) \, dt - \frac{\widehat{e_\star}}{2}\int_{1/2}^1 g(t) \, dt \right) \\
  &= (t\chi_{[0, 1/2]} + (1 - t)\chi_{[1/2, 1]})\left(\int_0^{1/2} g(t) \, dt - \int_{1/2}^1 g(t) \, dt \right)
\end{align*}
and thus $(\myoset{\star}{W} a^{\star} \myoset{\star}{R} + \myoset{\star}{s})\mid_{\Omega^{\bullet}(\left|\Delta^1\right|)} = s$.
\end{proof}

In general, the simplicial complex $\star_{i_0, \dots, i_k} \Delta^n$ has $k + 1$ top dimensional simplices $[e_{\star}, e_0, \dots, \widehat{e_{i_m}}, \dots, e_n]$ for $m = 0, \dots, k$.  For each top dimensional simplex in $\star_{i_0, \dots, i_k} \Delta^n$, there are barycentric coordinates which we would like to relate to barycentric coordinates on $\Delta^n$.  

Writing this down explicitly, a point $t_0e_0 + \dots + t_ne_n$ in $\Delta^n$ is in the $m$-th top dimensional simplex of $\star_{i_0, \dots, i_k}\Delta^n$ if it is a convex combination
\begin{align*}
  t_\star'e_{\star} + t_0'e_0 + \dots + t_{i_m - 1}'e_{i_m - 1} + t_{i_m + 1}'e_{i_m + 1}+ \dots + t_n'e_n
\end{align*}
where $e_{\star} = (e_{i_0} + \dots + e_{i_k})/(k + 1)$.  Or equivalently, if
\begin{align*}
  t_{i_m} \leq 1/(k + 1) \quad \text{and} \quad t_i \geq t_{i_m} \text{ for all $i \in I \setminus \{i_m\}$}.
\end{align*}
The change of barycentric coordinates is given by
\begin{align*}
  t_{\star}' = (k + 1)t_{i_m} \quad \text{and} \quad t_i' = t_i - t_{i_m} \quad \text{and} \quad t_j' = t_j 
\end{align*}
for $i \in I \setminus \{i_m\}$ and for $j \not \in I$.

\begin{myprop}
  On $[e_{\star}, e_0, \dots, \widehat{e_{i_m}}, \dots, e_n]$, the $m$-th top dimensional simplex of $\star_{i_0, \dots, i_k}\Delta^n$, the Whitney forms are denoted $\overline{\omega}_{j_0, \dots, j_l}'$ and $\overline{\omega}_{\star, j_0, \dots, j_l}'$, for all $J \not \ni i_m$.  These are related to the Whitney forms on $\Delta^n$ by the formulas
  \begin{align*}
    \overline{\omega}_{j_0, \dots, j_l}' = \overline{\omega}_{j_0, \dots, j_l} - \sum_{j_i \in I} (-1)^i\overline{\omega}_{i_m, j_0, \dots, \widehat{j_i}, \dots, j_l}
  \end{align*}
  and
  \begin{align*}
    \overline{\omega}_{\star, j_0, \dots, j_l}' = (k + 1)\overline{\omega}_{i_m, j_0, \dots, j_l}
  \end{align*}
\end{myprop}
\begin{proof}
On $[e_{\star}, e_0, \dots, \widehat{e_{i_m}}, \dots, e_n]$, we compute that
\begin{align*}
  \overline{\omega}_{j_0, \dots, j_l}' &=  \overline{\omega}_{j_0, \dots, j_l} -  l!\sum_{j_i \in I} (-1)^i\left(t_{i_m} dt_{j_0} \dots \widehat{dt_{j_i}} \dots dt_{j_l} -  dt_{i_m} \omega_{j_0, \dots, \widehat{j_i}, \dots, j_l} \right)\\
  &= \overline{\omega}_{j_0, \dots, j_l} - \sum_{j_i \in I} (-1)^i\overline{\omega}_{i_m, j_0, \dots, \widehat{j_i}, \dots, j_l}.
\end{align*}
Also, on $[e_{\star}, e_0, \dots, \widehat{e_{i_m}}, \dots, e_n]$, the identities
\begin{align*}
  dt_{j_0}' \dots dt_{j_l}' &= dt_{j_0} \dots dt_{j_l} -  dt_{i_m} \sum_{j_i \in I} (-1)^idt_{j_0} \dots \widehat{dt_{j_i}} \dots  dt_{j_l}\\
  dt_\star' \omega_{j_0, \dots, j_l}' &= (k + 1) dt_{i_m} \omega_{j_0, \dots, j_l} - (k + 1)t_{i_m} dt_{i_m} \sum_{j_i \in I} (-1)^i dt_{j_0} \dots \widehat{dt_{j_i}} \dots dt_{j_l}
\end{align*}
imply the formula
\begin{align*}
  \overline{\omega}_{\star, j_0, \dots, j_l}' &= (l + 1)!\left(t_{\star}' dt_{j_0}' \dots dt_{j_l}' - dt_{\star}'\omega_{j_0, \dots, j_l}'\right)\\
                                              &= (l + 1)!(k + 1)\left( t_{i_m} dt_{j_0} \dots dt_{j_l} -  dt_{i_m} \omega_{j_0, \dots, j_l}\right) \\
  &= (k + 1)\overline{\omega}_{i_m, j_0, \dots, j_l}
\end{align*}
\end{proof}

On $\star_{i_0, \dots, i_k}\Delta^n$, there is a deformation retraction
\begin{align*}
  \myDR{(C^{\bullet}(\star_{i_0, \dots, i_k}\Delta^n), d)}{(\Omega^{\bullet}(\left|\star_{i_0, \dots, i_k}\Delta^n\right|), d)}{\accentset{\star}{W}}{\accentset{\star}{R}}{\accentset{\star}{s}}
\end{align*}

\begin{mythm}\label{thm:dup_stell_subdiv}
Composing with the stellar subdivision deformation retraction gives a new deformation retraction
\begin{align*}
  \myDR{(C^{\bullet}(\Delta^n), d)}{(\Omega^{\bullet}(\left|\star_{i_0, \dots, i_k}\Delta^n\right|), d)}{\accentset{\star}{W} i^{\star}}{p^{\star} \accentset{\star}{R}}{\accentset{\star}{s} + \accentset{\star}{W} a^{\star} \accentset{\star}{R}}
\end{align*}
Comparing this to the deformation retraction,
\begin{align*}
  \myDR{(C^{\bullet}(\Delta^n), d)}{(\Omega^{\bullet}(\left|\Delta^n\right|), d)}{W}{R}{s}
\end{align*}
we find that $p^{\star} \myoset{\star}{R}|_{\Omega^{\bullet}(|\Delta^n|)} = R$, $\myoset{\star}{W} i^{\star} = W$, and $(\myoset{\star}{s} + \myoset{\star}{W} a^{\star} \myoset{\star}{R})|_{\Omega^{\bullet}(|\Delta^n|)} - s$ is closed.
\end{mythm}
\begin{proof}
We have
\begin{align*}
  p^{\star} \myoset{\star}{R} &= p^{\star}\left[\sum_{\substack{j_0 < \dots < j_l\\ J \not \supset I}} \reallywidehat{[j_0, \dots, j_l]} \int_{[j_0, \dots, j_l]} + \sum_{\substack{j_0 < \dots < j_l\\ J \not \supset I}} \reallywidehat{[\star, j_0, \dots, j_l]} \int_{[\star, j_0, \dots, j_l]}\right]\\
                               &= \sum_{\substack{j_0 < \dots < j_l\\ J \not \supset I}} \reallywidehat{[j_0, \dots, j_l]} \int_{[j_0, \dots, j_l]} + \sum_{\substack{j_0 < \dots < j_l\\ J \supset I}} \reallywidehat{[j_0, \dots, j_l]} \int_{[j_0, \dots, j_l]}\\
  &= R
\end{align*}

Secondly, for $J \supset I$, we have on $[e_{\star}, e_0, \dots, \widehat{e_{i_m}}, \dots, e_n]$ for $i_m = j_b$
\begin{align*}
  \myoset{\star} W i^{\star}\reallywidehat{[e_{j_0}, \dots, e_{j_l}]} &= \frac{1}{k + 1}\sum_{j_i \in I} (-1)^i\overline{\omega}_{\star, j_0, \dots, \widehat{j_i}, \dots, j_l}' = \frac{(-1)^b}{k + 1} \overline{\omega}_{\star, j_0, \dots, \widehat{j_b}, \dots, j_l}' \\                            
                                                                       &= (-1)^b\overline{\omega}_{j_b, j_0, \dots, \widehat{j_b}, \dots, j_l} = \overline{\omega}_{j_0, \dots, j_l} =  W\reallywidehat{[e_{j_0}, \dots, e_{j_l}]}
\end{align*}
For $J \not \supset I$, we have on $[e_{\star}, e_0, \dots, \widehat{e_{i_m}}, \dots, e_n]$ if $i_m \in I \setminus J$
\begin{align*}
  \myoset{\star} W i^{\star}\reallywidehat{[e_{j_0}, \dots, e_{j_l}]} &= \overline{\omega}_{j_0, \dots, j_l}' + \frac{1}{k + 1} \sum_{j_i \in I} (-1)^i \overline{\omega}_{\star, j_0, \dots, \widehat{j_i}, \dots, j_l}'\\
                                                                       &= \overline{\omega}_{j_0, \dots, j_l} -  \sum_{j_i \in I} (-1)^i \overline{\omega}_{i_m, j_0, \dots, \widehat{j_i}, \dots, j_l} + \sum_{j_i \in I}(-1)^i \overline{\omega}_{i_m, j_0, \dots, \widehat{j_i}, \dots, j_l}\\
  &= \overline{\omega}_{j_0, \dots, j_l}  = W\reallywidehat{[e_{j_0}, \dots, e_{j_l}]}
\end{align*}
For $J \not \supset I$, we have on $[e_{\star}, e_0, \dots, \widehat{e_{i_m}}, \dots, e_n]$ if $i_m = j_b$.
\begin{align*}
  \myoset{\star} W i^{\star}\reallywidehat{[e_{j_0}, \dots, e_{j_l}]} &= \overline{\omega}_{j_0, \dots, j_l}' + \frac{1}{k + 1} \sum_{j_i \in I} (-1)^i \overline{\omega}_{\star, j_0, \dots, \widehat{j_i}, \dots, j_l}' \\
                                                                       &= \frac{(-1)^b }{k + 1} \overline{\omega}_{\star, j_0, \dots, \widehat{j_b}, \dots, j_l}'\\
  &= \overline{\omega}_{j_0, \dots, j_l} =  W\reallywidehat{[e_{j_0}, \dots, e_{j_l}]}
\end{align*}
Solely based on the identities $p^{\star} \myoset{\star}{R}|_{\Omega^{\bullet}(|\Delta^n|)} = R$ and $\myoset{\star}{W} i^{\star} = W$, it follows that $(\myoset{\star}{s} + \myoset{\star}{W} a^{\star} \myoset{\star}{R})|_{\Omega^{\bullet}(|\Delta^n|)} - s$ is closed.
\end{proof}
In the rest of the section, we investigate the structure of $\myoset{\star} W a^{\star} \myoset{\star} R + \myoset{\star}s - s$.  We begin by turning to the Dupont homotopy $\myoset{\star}{s}$.  We find on each top dimensional simplex $[e_{\star}, e_0, \dots, \widehat{e_{i_m}}, \dots, e_n]$ that
\begin{align*}
  \phi^{\star'}\left(s, \sum_jt_je_j\right) &:= (1 - s)\left[t_{\star}'e_{\star} + \sum_{j \neq i_m} t_j' e_j\right] + se_\star \\
  &= \sum_j (1 - s)t_je_j + se_\star = \phi^{\star}\left(s, \sum_jt_je_j\right)
\end{align*}
where we have used the fact that $e_{\star} = \frac{1}{k + 1}\sum_{i \in I} se_i$.  Also on each top dimensional simplex $[e_{\star}, e_0, \dots, \widehat{e_{i_m}}, \dots, e_n]$, for $\alpha \neq i_m$ we have that
\begin{align*}
  \phi^{\alpha'}\left(s, \sum_jt_je_j\right) &:= (1 - s)\left[t_{\star}'e_{\star} + \sum_{j \neq i_m} t_j' e_j\right] + se_\alpha \\
  &= \sum_j (1 - s)t_je_j + se_\alpha = \phi^{\alpha}\left(s, \sum_jt_je_j\right) .
\end{align*}

The Dupont homotopy is given by
\begin{align*}
  \myoset{\star} s = -\sum_{l = 0}^{n - 1} \sum_{\substack{j_0 < \dots < j_l \\ J \not \supset I}} \overline{\omega}_{j_0, \dots, j_l}' h^{j_l'}\dots h^{j_0'} -\sum_{l = -1}^{n - 2} \sum_{\substack{j_0 < \dots < j_l \\ J \not \supset I}} \overline{\omega}_{\star, j_0, \dots, j_l}' h^{j_l'}\dots h^{j_0'}h^{\star'} 
\end{align*}
For $l = -1$ we are simply defining the inner sum to be equal to $\overline{\omega}_{\star} h^{\star}$.  
On $[e_{\star}, e_0, \dots, \widehat{e_{i_m}}, \dots, e_n]$ the Dupont homotopy restricts to
\begin{align*}
  \myoset{\star} s &= -\sum_{l = 0}^{n - 1} \sum_{\substack{j_0 < \dots < j_l \\ i_m \not \in J}} \overline{\omega}_{j_0, \dots, j_l}' h^{j_l}\dots h^{j_0} - \sum_{l = -1}^{n - 2} \sum_{\substack{j_0 < \dots < j_l \\ i_m \not \in J}} \overline{\omega}_{\star, j_0, \dots, j_l}' h^{j_l}\dots h^{j_0} h^{\star} \\
                    &= -\sum_{l = 0}^{n - 1} \sum_{\substack{j_0 < \dots < j_l \\ i_m \not \in J}} \overline{\omega}_{j_0, \dots, j_l} h^{j_l} \dots h^{j_0} \\
  &+ \sum_{l = -1}^{n - 2} \sum_{\substack{j_0 < \dots < j_l \\ i_m \not \in J}} \overline{\omega}_{i_m, j_0, \dots j_l} h^{j_l} \dots h^{j_0}\left( \sum_{i_{\alpha} \in I \setminus (\{i_m\} \cup J) } h^{i_{\alpha}} - (k + 1)h^{\star}\right)
\end{align*}

\begin{mythm}
  \begin{align*}
    \myoset{\star} W a^{\star} \myoset{\star} R + \myoset{\star}s - s = \sum_{i_m \in I}\sum_{l = -1}^{n - 2} \sum_{\substack{j_0 < \dots < j_l \\ i_m \not \in J}} \chi_{i_m}\overline{\omega}_{i_m, j_0, \dots, j_l} T_{j_0, \dots, j_l}
  \end{align*}
  where
  \begin{align*}
    T_{j_0, \dots, j_l} &= (-1)^{l + 1}\left((1 - \varepsilon^{\star}) \sum_{i_{\alpha} \in I \setminus J} h^{i_{\alpha}} - (k + 1)h^{\star} \right) h^{j_l} \dots h^{j_0}
  \end{align*}
  Here $\varepsilon^{\star}$ is evaluation at $e_{\star}$, the barycenter of $[i_0, \dots, i_k]$ and $\chi_{i_m}$ is the characteristic function for the subset $[e_{\star}, e_0, \dots, \widehat{e_{i_m}}, \dots, e_n]$.
\end{mythm}
\begin{proof}

We have
\begin{align*}
  \myoset{\star}{W} a^{\star} \myoset{\star}{R} &= -\frac{1}{k + 1}\sum_{l = 0}^{n - 1} \sum_{\substack{j_0 < \dots < j_l\\ J \not \supset I}}  \sum_{j_i \in I} (-1)^i \overline{\omega}'_{\star, j_0, \dots, \widehat{j_i}, \dots, j_l} \int_{[\star, j_0, \dots, j_l]}
\end{align*}
On $[e_{\star}, e_0, \dots, \widehat{e_{i_m}}, \dots, e_n]$, this reduces to
\begin{align*}
  \myoset{\star}{W} a^{\star} \myoset{\star}{R} &= -\frac{1}{k + 1}\sum_{l = 0}^{n - 1} \sum_{\substack{j_0 < \dots < j_l\\ i_m \not \in J}}  \sum_{j_i \in I} (-1)^i \overline{\omega}'_{\star, j_0, \dots, \widehat{j_i}, \dots, j_l} \int_{[\star, j_0, \dots, j_l]}\\
                                                  &-\frac{1}{k + 1}\sum_{l = 0}^{n - 1} \sum_{\substack{j_0 < \dots < j_l\\ i_m = j_b, J \not \supset I}}  (-1)^b \overline{\omega}'_{\star, j_0, \dots, \widehat{j_b}, \dots, j_l} \int_{[\star, j_0, \dots, j_l]}\\
  &= -\sum_{l = -1}^{n - 2} \sum_{\substack{j_0 < \dots < j_l\\ i_m \not \in J}}   \overline{\omega}_{i_m, j_0, \dots, j_l} \sum_{\alpha \in I \setminus (\{i_m\} \cup J)} \int_{[\star, i_{\alpha}, j_0, \dots, j_l]}\\
  &-\sum_{l = -1}^{n - 2} \sum_{\substack{j_0 < \dots < j_l\\ i_m \not \in J}}   \overline{\omega}_{i_m, j_0, \dots, j_l} \int_{[\star, i_m, j_0, \dots, j_l]}
\end{align*}
Therefore
\begin{align*}
  \myoset{\star}{W} a^{\star} \myoset{\star}{R} = -\sum_{l = -1}^{n - 2} \sum_{\substack{j_0 < \dots < j_l\\ i_m \not \in J}}   \overline{\omega}_{i_m, j_0, \dots, j_l} \sum_{\alpha \in I \setminus J} \int_{[\star, i_{\alpha}, j_0, \dots, j_l]}
\end{align*}

Collecting the results gives
\begin{align*}
  T_{j_0, \dots, j_l} &= -\sum_{\alpha \in I \setminus J} \int_{[\star, i_{\alpha}, j_0, \dots, j_l]} \\
                      &+ h^{j_l} \dots h^{j_0}\left( \sum_{i_{\alpha} \in I \setminus (\{i_m\} \cup J) } h^{i_{\alpha}} - (k + 1)h^{\star} + h^{i_m} \right).
\end{align*}
Bringing $h^{j_l} \dots h^{j_0}$ to the right gives the desired result.
\begin{align*}
  T_{j_0, \dots, j_l} = (-1)^{l + 1}\left(-\sum_{i_{\alpha} \in I \setminus J} \varepsilon^{\star}h^{i_{\alpha}} + \sum_{i_{\alpha} \in I \setminus J} h^{i_{\alpha}} - (k + 1)h^{\star} \right) h^{j_l} \dots h^{j_0}
\end{align*}
\end{proof}

As a sanity check, we can verify that $T_{\emptyset}(\omega) = 0$ for all $\omega$, when $n = 1$, $k = 1$.  That is,
\begin{align*}
  (1 - \varepsilon^{\star}) \sum_{i_{\alpha} \in I} h^{i_{\alpha}}(\omega) = (k + 1)h^{\star}(\omega),
\end{align*}
which for $n = 1$ becomes the formula
\begin{align*}
   (1 - \varepsilon^{\star})(h^0 + h^1)(g(t)\, dt) &= (1 - \varepsilon^{\star}) \left[\int_t^0 g(u) \, du + \int_t^1 g(u) \, du \right] \\
                                                &= 2 \int_t^{1/2} g(u) \, du = 2 h^{\star}(g(t) \, dt)
\end{align*}
as expected.

For general $n$, it is not true that $\myoset{\star} W a^{\star} \myoset{\star} R + \myoset{\star}s - s = 0$.  However, $(\myoset{\star} W a^{\star} \myoset{\star} R)|_{\Omega^{\bullet}(|\Delta^n|)} + \myoset{\star}s - s$ is exact.  That is, $(\myoset{\star}{s} + \myoset{\star}{W} a^{\star} \myoset{\star}{R})|_{\Omega^{\bullet}(|\Delta^n|)}  - s = dQ - Qd$ for some map $Q$.  This is because $\Delta^n$ is contractible, so the homomorphism complex $\Hom^\bullet(\Omega(|\Delta^n|), \Omega(|\Delta^n|))$ defined by
\begin{align*}
  \Hom^i(\Omega(|\Delta^n|), \Omega(|\Delta^n|)) = \oplus_{j = 0}^n \Hom(\Omega^j(|\Delta^n|), \Omega^{j + i}(|\Delta^n|))
\end{align*}
has cohomology
\begin{align*}
  H^i(\Hom(\Omega(|\Delta^n|), \Omega(|\Delta^n|))) \cong
  \begin{cases}
    \Bbbk \text{ if $i = 0$} \\
    0 \text{ if $i \neq 0$}
  \end{cases}
\end{align*}
Therefore a closed degree $-1$ linear endomorphism of $\Omega^{\bullet}(|\Delta^n|)$ must be exact.
The best we can do is give an explicit formula for the homotopy $Q$.  We leave this as an exercise for the reader.

\section{Compatibility of the Cubical DHF and Cubical Stellar Subdivision}
\label{sec:comp_cub_stell}
The Dupont homotopy on $\star_{1, \dots, k} \square^n$ is a deformation retraction
\begin{align*}
  \myDR{(C^{\bullet}(\star_{1, \dots, k} \square^n), d)}{(\Omega^{\bullet}(\star_{1, \dots, k} \square^n), d)}{\myoset{\star}{\mathbf{I}}}{\myoset{\star}{\mathbf{W}}}{\myoset{\star}{\mathbf{s}}}
\end{align*}
where $\myoset{\star}{\mathbf{I}} = (\myoset{\star} I)^{\otimes k} \otimes I^{\otimes (n - k)} = \myoset{\star}{\mathbf{I}}_k \otimes \mathbf{I}_{n - k}$, $\myoset{\star}{\mathbf{W}} = (\myoset{\star} W)^{\otimes k} \otimes W^{\otimes (n - k)} = \myoset{\star}{\mathbf{W}}_k \otimes \mathbf{W}_{n - k}$ and
\begin{align*}
  \myoset{\star}{\mathbf{s}} = \frac{1}{n!}\sum_{\mathbf{\epsilon}} C_{|\epsilon|, n} \boldsymbol{\psi}_{\epsilon}
\end{align*}
where $C_{|\epsilon|, n} = |\epsilon|!(n - 1 - |\epsilon|)!$ and the outer sum is over $\epsilon_k = 0, 1$ and where
\begin{align*}
  \boldsymbol{\psi}_{\epsilon} &= \sum_{j = 1}^k \underbrace{(\myoset{\star}W \myoset{\star}I)^{\epsilon_1} \otimes \dots}_{j - 1}  \otimes \, \myoset{\star} s \,\otimes \dots \otimes (\myoset{\star}W \myoset{\star} I)^{\epsilon_{k - 1}} \otimes (WI)^{\epsilon_k} \otimes \dots \otimes (WI)^{\epsilon_{n - 1}}\\
  &+\sum_{j = 1}^{n - k} (\myoset{\star}W \myoset{\star}I)^{\epsilon_1} \otimes \dots \otimes (\myoset{\star}W \myoset{\star}I)^{\epsilon_k}  \otimes \underbrace{(WI)^{\epsilon_{k + 1}} \otimes \dots}_{j - 1} \otimes \, s \, \otimes \dots \otimes (WI)^{\epsilon_{n - 1}}
\end{align*}

\begin{mythm}
  Composing with the stellar subdivision deformation retraction gives a new deformation retraction
  \begin{align*}
    \myDR{(C^{\bullet}(\square^n), d)}{(\Omega^{\bullet}(\star_{1, \dots, k}\square^n), d)}{\myoset{\star}{\mathbf{W}} \mathbf{i}^{\star}}{\mathbf{p}^{\star}\myoset{\star}{\mathbf{I}}}{\myoset{\star}{\mathbf{s}} + \myoset{\star}{\mathbf{W}} \mathbf{a}^{\star} \myoset{\star}{\mathbf{I}}}
  \end{align*}
  Comparing this to the deformation retraction
  \begin{align*}
    \myDR{(C^{\bullet}(\square^n), d)}{(\Omega^{\bullet}(\square^n), d)}{\mathbf{W}}{\mathbf{I}}{\mathbf{s}}
  \end{align*}
  we find that $\mathbf{p}^{\star}\myoset{\star}{\mathbf{I}}|_{\Omega^{\bullet}(\square^n)} = \mathbf{I}$, $\myoset{\star}{\mathbf{W}} \mathbf{i}^{\star} = \mathbf{W}$, and $(\myoset{\star}{\mathbf{W}} \mathbf{a}^{\star} \myoset{\star}{\mathbf{I}} + \myoset{\star}{\mathbf{s}})|_{\Omega^{\bullet}(\square^n)} - \mathbf{s}$ is closed.
\end{mythm}
\begin{proof}
  The identities $\mathbf{p}^{\star}\myoset{\star}{\mathbf{I}} = \mathbf{I}$ and $\myoset{\star}{\mathbf{W}} \mathbf{i}^{\star} = \mathbf{W}$ follow from the identities $p^{\star} \myoset{\star}I = I$ and $\myoset{\star}W i^{\star} = W$ on $\Delta^1$ that were proved in Theorem \ref{thm:comp-dhf-stell-1}.

  We shall try to write down the homotopy $\mathbf{Q}$ explicitly only in the more symmetric case $k = n$ for which the formulas are easier to handle.  For $k = n$,
  \begin{align*}
    \myoset{\star}{\mathbf{W}} \mathbf{a}^{\star} \myoset{\star}{\mathbf{I}} + \myoset{\star}{\mathbf{s}} - \mathbf{s} &= \sum_{\sigma \in S_n} \tau_{\sigma} \circ \left(\sum_{j = 1}^n (\myoset{\star}W \myoset{\star}I)^{\otimes (j - 1)} \otimes \myoset{\star}W a^{\star} \myoset{\star}I \otimes (WI)^{\otimes(n - j)} \right) \circ \tau_{\sigma}^{-1}\\
                                                                                                                          &+ \sum_{\sigma \in S_n} \tau_{\sigma} \circ \left(\sum_{j = 1}^n 1^{\otimes (j - 1)} \otimes \myoset{\star}s \otimes (\myoset{\star}W \myoset{\star}I)^{\otimes(n - j)}  \right) \circ \tau_{\sigma}^{-1}\\
    &- \sum_{\sigma \in S_n} \tau_{\sigma} \circ \left(\sum_{j = 1}^n 1^{\otimes (j - 1)} \otimes s \otimes (WI)^{\otimes(n - j)} \right) \circ \tau_{\sigma}^{-1}
  \end{align*}
  And thus
  \begin{align*}
    \myoset{\star}{\mathbf{W}} \mathbf{a}^{\star} \myoset{\star}{\mathbf{I}} + \myoset{\star}{\mathbf{s}} - \mathbf{s} &= \sum_{\sigma \in S_n} \tau_{\sigma} \circ \left(\sum_{j = 1}^n (\myoset{\star}W \myoset{\star}I - 1) ^{\otimes (j - 1)} \otimes s \otimes (WI)^{\otimes(n - j)} \right) \circ \tau_{\sigma}^{-1} \\
    &+ \sum_{\sigma \in S_n} \tau_{\sigma} \circ \left(\sum_{j = 1}^n (\myoset{\star}W \myoset{\star}I) ^{\otimes (j - 1)} \otimes \myoset{\star}s \otimes (1 - WI)^{\otimes(n - j)} \right) \circ \tau_{\sigma}^{-1}
  \end{align*}
\end{proof}

\section{Elementary Expansion and Collapse}
\label{sec:el_coll}
Let $Y$ be a simplicial complex containing a $k$-simplex $\sigma$ and a $(k - 1)$-simplex $\sigma'$ such that $\sigma$ is the only $k$-simplex whose boundary contains $\sigma'$.  Let $X \subset Y$ be the subcomplex obtained by removing the pair $\sigma$, $\sigma'$ from $Y$.  Then one calls $X$ an elementary collapse of $Y$ and $Y$ an elementary expansion of $X$.

To write down simplicial chains, we choose an orientation for each simplex.  Suppose that $\partial \sigma = \sum \varepsilon_{\tau} \tau$ where $\varepsilon_{\tau} = \pm 1$.  There is a natural projection $p_{\downarrow}: C_{\bullet}(Y) \to C_{\bullet}(X)$ sending $\sigma' \mapsto \sigma' - \varepsilon_{\sigma'}\partial\sigma$ and $\sigma \mapsto 0$.

\begin{myprop}\label{prop:ell-coll-chains}
  There is an elementary collapse deformation retraction
  \begin{align*}
    \myDR{C_{\bullet}(X)}{C_{\bullet}(Y)}{i_{\downarrow}}{p_{\downarrow}}{a_{\downarrow}}
  \end{align*}
  where $i_{\downarrow}$ is the natural inclusion, $p_{\downarrow}$ is as above, and $a_{\downarrow}(\sigma') = \varepsilon_{\sigma'}\sigma$ and $a_{\downarrow}(\tau) = 0$ for $\tau \neq \sigma'$.  
\end{myprop}
\begin{proof}
  We verify that $p_{\downarrow}$ is a chain map by computing $\partial p_{\downarrow}(\sigma') = \partial(\sigma' - \varepsilon_{\sigma'}\partial \sigma) = p_{\downarrow}\partial (\sigma')$ and
  \begin{align*}
    \varepsilon_{\sigma'}p_{\downarrow} \partial(\sigma) &= p_{\downarrow} (\sigma') + p_{\downarrow}(\varepsilon_{\sigma'}\partial(\sigma) - \sigma')  \\
    &= \sigma' - \varepsilon_{\sigma'}\partial(\sigma) + \varepsilon_{\sigma'}\partial(\sigma) - \sigma' = 0 = \partial p_{\downarrow}(\sigma).
  \end{align*}  

  Lastly, we verify that $\partial a_{\downarrow} + a_{\downarrow} \partial = 1 - i_{\downarrow}p_{\downarrow}$.  We have
\begin{align*}
  \partial a_{\downarrow}(\sigma) + a_{\downarrow}\partial (\sigma) = a_{\downarrow}\partial (\sigma) = \sigma = \sigma - i_{\downarrow}p_{\downarrow}(\sigma)
\end{align*}
and
\begin{align*}
  \partial a_{\downarrow}(\sigma') + a_{\downarrow}\partial (\sigma') = \partial a_{\downarrow} (\sigma') = \varepsilon_{\sigma'}\partial \sigma = \sigma' - (\sigma' - \varepsilon_{\sigma'}\partial \sigma) = \sigma' - i_{\downarrow}p_{\downarrow}(\sigma').
\end{align*}

\end{proof}

Let 
\begin{align*}
  i^{\downarrow}(\widehat{\tau}) &=
  \begin{cases}
    \widehat{\tau} - \varepsilon_{\sigma'}\varepsilon_{\tau}\widehat{\sigma'} & \text{if $\tau$ is in $\partial \sigma - \sigma'$}\\
    \widehat{\tau} & \text{otherwise}
  \end{cases}\\
  p^{\downarrow}(\widehat{\tau}) &=
                      \begin{cases}
                        \widehat{\tau} & \text{if $\tau$ is in $X$}\\
                        0 & \text{if $\tau = \sigma, \sigma'$}
                      \end{cases}
\end{align*}
and lastly,
\begin{align*}
  a^{\downarrow}(\widehat{\tau}) &=
                      \begin{cases}
                        \varepsilon_{\sigma'}\widehat{\sigma'} & \text{if $\tau = \sigma$}\\
                        0 & \text{otherwise}
                      \end{cases}
\end{align*}
\begin{myprop}
  There is an elementary collapse deformation retraction on cochains
  \begin{align*}
    \myDR{C^{\bullet}(X)}{C^{\bullet}(Y)}{i^{\downarrow}}{p^{\downarrow}}{a^{\downarrow}}
  \end{align*}
  where $i^{\downarrow}$, $p^{\downarrow}$ and $a^{\downarrow}$ are as defined above.  This is dual to the deformation retraction defined in Proposition \ref{prop:ell-coll-chains}.
\end{myprop}

\section{Stellar Subdivision from Elementary Expansions and Collapses}
\label{sec:stell_subdiv2}
A stellar subdivision of the $n$-simplex can be constructed as a sequence of elementary expansions and elementary collapses.  For example for the $1$-simplex, there are two such sequences:
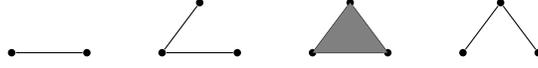
\begin{center}
  \begin{tikzpicture}
    \tikzstyle{vertex} = [circle, fill=black, inner sep=1pt]
    \node[vertex] (A) at (0, 0) {};
    \node[vertex] (B) at (1, 0) {};
    \draw (A) -- (B);
    
    \node[vertex] (A) at (2, 0) {};
    \node[vertex] (B) at (3, 0) {};
    \node[vertex] (C) at (2.5, .6666) {};
    \draw (A) -- (B);
    \draw (A) -- (C);

    \node[vertex] (A) at (4, 0) {};
    \node[vertex] (B) at (5, 0) {};
    \node[vertex] (C) at (4.5, .6666) {};
    \draw (A) -- (C) -- (B) -- (A);
    \fill[gray] (A.center) -- (C.center) -- (B.center) -- (A.center);

    \node[vertex] (A) at (6, 0) {};
    \node[vertex] (B) at (7, 0) {};
    \node[vertex] (C) at (6.5, .6666) {};
    \draw (A) -- (C) -- (B);
  \end{tikzpicture}
\end{center}
and
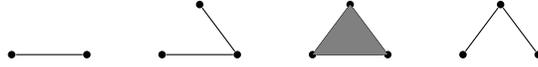
\begin{center}
  \begin{tikzpicture}
    \tikzstyle{vertex} = [circle, fill=black, inner sep=1pt]
    \node[vertex] (A) at (0, 0) {};
    \node[vertex] (B) at (1, 0) {};
    \draw (A) -- (B);
    
    \node[vertex] (A) at (2, 0) {};
    \node[vertex] (B) at (3, 0) {};
    \node[vertex] (C) at (2.5, .6666) {};
    \draw (A) -- (B);
    \draw (B) -- (C);

    \node[vertex] (A) at (4, 0) {};
    \node[vertex] (B) at (5, 0) {};
    \node[vertex] (C) at (4.5, .6666) {};
    \draw (A) -- (C) -- (B) -- (A);
    \fill[gray] (A.center) -- (C.center) -- (B.center) -- (A.center);

    \node[vertex] (A) at (6, 0) {};
    \node[vertex] (B) at (7, 0) {};
    \node[vertex] (C) at (6.5, .6666) {};
    \draw (A) -- (C) -- (B);
  \end{tikzpicture}
\end{center}
We start by identifying $\Delta^1 = [e_0, e_1]$ in $\Delta^2$ as $\Delta^2 = [e_0, e_\star, e_1]$.  All simplices that we consider in the sequence of elementary expansions and collapses shall be given the induced orientation from $\Delta^2$.  We shall indicate the simplex $\sigma$ in the subscript of the maps in the elementary collapse deformation retractions, and use primes to indicate the face $\sigma' \subset \sigma$.

We consider the first sequence of elementary expansions.  The first inclusion map is
\begin{align*}
  i^{[e_0, e'_\star]}(\widehat{\tau}) =
  \begin{cases}
    \widehat{e_0} + \widehat{e_\star} & \text{if $\tau = e_0$}\\
    \widehat{\tau} & \text{otherwise}
  \end{cases}
\end{align*}
and the second inclusion map is
\begin{align*}
  i^{[e_0, e'_{\star}, e'_1]}(\widehat{\tau}) =
  \begin{cases}
    \widehat{[e_0, e_{\star}]} - \widehat{[e_{\star}, e_1]} & \text{if $\tau = [e_0, e_{\star}]$}\\
    \widehat{[e_0, e_1]} + \widehat{[e_{\star}, e_1]} & \text{if $\tau = [e_0, e_1]$}\\
    \widehat{\tau} & \text{otherwise}
  \end{cases}
\end{align*}
and their composition is the inclusion map
\begin{align*}
  i_1^{\downarrow}(\widehat{\tau}) = i^{[e_0, e'_{\star}, e'_1]}i^{[e_0, e'_{\star}]}(\widehat{\tau}) =
  \begin{cases}
    \widehat{[e_0, e_1]} + \widehat{[e_{\star}, e_1]} & \text{if $\tau = [e_0, e_1]$}\\
    \widehat{e_0} + \widehat{e_{\star}} & \text{if $\tau = e_0$}\\
    \widehat{\tau} & \text{otherwise}
  \end{cases}
\end{align*}
the projection map 
\begin{align*}
  p_1^{\downarrow}(\widehat{\tau}) = p^{[e_0, e'_{\star}]}p^{[e_0, e'_{\star}, e'_1]}(\widehat{\tau}) =
  \begin{cases}
    \widehat{\tau} &\text{if $\tau \in \Delta^1$}\\
    0 & \text{otherwise}
  \end{cases}
\end{align*}
and the homotopy
\begin{align*}
  a_1^{\downarrow}(\widehat{\tau}) = (a^{[e_0, e'_{\star}, e'_1]} + i^{[e_0, e'_{\star}, e'_1]}a^{[e_0, e'_{\star}]}p^{[e_0, e'_{\star}, e'_1]})(\widehat{\tau}) =
  \begin{cases}
    \widehat{[e_{\star}, e_1]} & \text{if $\tau = [e_0, e_{\star}, e_1]$}\\
    \widehat{e_{\star}} & \text{if $\tau = [e_0, e_{\star}]$}\\
    0 & \text{otherwise}
  \end{cases}
\end{align*}
We then apply the elementary collapse to get
\begin{align*}
  \myoset{\star}{\iota}_1(\widehat{\tau}) = p^{[e'_0, e_{\star}, e'_1]}i_1^{\downarrow}(\widehat{\tau}) =
  \begin{cases}
    \widehat{[e_{\star}, e_1]} & \text{if $\tau = [e_0, e_1]$}\\
    \widehat{e_0} + \widehat{e_{\star}} & \text{if $\tau = e_0$}\\
    \widehat{\tau} & \text{otherwise}
  \end{cases}
\end{align*}
and
\begin{align*}
  \myoset{\star}{p}_1(\widehat{\tau}) = p_1^{\downarrow} i^{[e'_0, e_{\star}, e'_1]}(\widehat{\tau}) =
  \begin{cases}
    \widehat{[e_0, e_1]} & \text{if $\tau = [e_0, e_{\star}]$}\\
    \widehat{[e_0, e_1]} & \text{if $\tau = [e_{\star}, e_1]$}\\
    0 & \text{if $\tau = e_{\star}$}\\
    \widehat{\tau} & \text{otherwise}
  \end{cases}
\end{align*}
and lastly
\begin{align*}
  \myoset{\star}{a}_1(\widehat{\tau}) = p^{[e'_0, e_{\star}, e'_1]}a_1^{\downarrow}i^{[e'_0, e_{\star}, e'_1]}(\widehat{\tau}) =
  \begin{cases}
    \widehat{e_{\star}} & \text{if $\tau = [e_0, e_{\star}]$}\\
    0 & \text{otherwise}.
  \end{cases}
\end{align*}

For the second sequence of elementary expansions:  The inclusion map is
\begin{align*}
  i_2^{\downarrow}(\widehat{\tau}) = i^{[e'_0, e'_{\star}, e_1]}i^{[e'_{\star}, e_1]}(\widehat{\tau}) =
  \begin{cases}
    \widehat{[e_0, e_1]} + \widehat{[e_0, e_{\star}]} & \text{if $\tau = [e_0, e_1]$}\\
    \widehat{e_1} + \widehat{e_{\star}} & \text{if $\tau = e_1$}\\
    \widehat{\tau} & \text{otherwise}
  \end{cases}
\end{align*}
the projection map is
\begin{align*}
  p_2^{\downarrow}(\widehat{\tau}) = p^{[e'_{\star}, e_1]}p^{[e'_0, e'_{\star}, e_1]}(\widehat{\tau}) =
  \begin{cases}
    \widehat{\tau} &\text{if $\tau \in \Delta^1$}\\
    0 & \text{otherwise}
  \end{cases}
\end{align*}
and the homotopy is
\begin{align*}
  a_2^{\downarrow}(\widehat{\tau}) = (a^{[e'_0, e'_{\star}, e_1]} + i^{[e'_0, e'_{\star}, e_1]}a^{[e'_{\star}, e_1]}p^{[e'_0, e'_{\star}, e_1]})(\widehat{\tau}) =
  \begin{cases}
    \widehat{[e_0, e_{\star}]} & \text{if $\tau = [e_0, e_{\star}, e_1]$}\\
    -\widehat{e_{\star}} & \text{if $\tau = [e_{\star}, e_1]$}\\
    0 & \text{otherwise}
  \end{cases}
\end{align*}
We then apply the elementary collapse to get
\begin{align*}
  \myoset{\star}{\iota}_2(\widehat{\tau}) = p^{[e'_0, e_{\star}, e'_1]}i_2^{\downarrow}(\widehat{\tau}) =
  \begin{cases}
    \widehat{[e_0, e_{\star}]} & \text{if $\tau = [e_0, e_1]$}\\
    \widehat{e_1} + \widehat{e_{\star}} & \text{if $\tau = e_1$}\\
    \widehat{\tau} & \text{otherwise}
  \end{cases}
\end{align*}
and
\begin{align*}
  \myoset{\star}{p}_2(\widehat{\tau}) = p_2^{\downarrow}i^{[e'_0, e_{\star}, e'_1]}(\widehat{\tau}) =
  \begin{cases}
    \widehat{[e_0, e_1]} & \text{if $\tau = [e_0, e_{\star}]$}\\
    \widehat{[e_0, e_1]} & \text{if $\tau = [e_{\star}, e_1]$}\\
    0 & \text{if $\tau = e_{\star}$}\\
    \widehat{\tau} & \text{otherwise}
  \end{cases}
\end{align*}
and lastly
\begin{align*}
  \myoset{\star}{a}_2(\widehat{\tau}) = p^{[e'_0, e_{\star}, e'_1]}a_2^{\downarrow}i^{[e'_0, e_{\star}, e'_1]} =
  \begin{cases}
    -\widehat{e_{\star}} & \text{if $\tau = [e_{\star}, e_1]$}\\
    0 & \text{otherwise}
  \end{cases}
\end{align*}

Let
\begin{align*}
  \myoset{\star}{\iota}(\widehat{\tau}) &= \frac{1}{2}(\myoset{\star}{\iota}_1 + \myoset{\star}{\iota}_2)(\widehat{\tau}) = 
  \begin{cases}
    \frac{1}{2}(\widehat{[e_0, e_{\star}]} + \widehat{[e_{\star}, e_1]}) & \text{if $\tau = [e_0, e_1]$}\\
    \frac{1}{2}\widehat{e_0} + \widehat{e_{\star}} & \text{if $\tau = e_0$}\\
    \frac{1}{2}\widehat{e_1} + \widehat{e_{\star}} & \text{if $\tau = e_1$}\\
    \widehat{\tau} & \text{otherwise},
  \end{cases}
\end{align*}
let $\myoset{\star}{p} = \myoset{\star}{p}_1 = \myoset{\star}{p}_2$ and let
\begin{align*}
  \myoset{\star}{a}(\widehat{\tau}) = \frac{1}{2}(\myoset{\star}{a}_1 + \myoset{\star}{a}_2)(\widehat{\tau}) =
  \begin{cases}
    \frac{1}{2}\widehat{e_1} & \text{if $\tau = [e_0, e_1]$}\\
    -\frac{1}{2}\widehat{e_1} & \text{if $\tau = [e_1, e_2]$}\\
    0 & \text{otherwise}
  \end{cases}
\end{align*}
Then we have $d\myoset{\star}{a} + \myoset{\star}{a}d = 1 - \myoset{\star}{\iota}\myoset{\star}{p}$, and remarkably we still have the property $\myoset{\star}{p}\myoset{\star}{\iota} = 1$, despite only composing a zigzag of deformations retractions.  More remarkably, this deformation retraction is equal to the stellar subdivision deformation retraction on $\Delta^1$ that we defined at the beginning of Section \ref{sec:stellar_subdiv_form}.  That is, $\myoset{\star}{\iota} = i^{\star}$, $\myoset{\star}{p} = p^{\star}$ and $\myoset{\star}{h} = h^{\star}$.

More generally, for $\star \Delta^n$, the stellar subdivision of the $n$-simplex, we shall give $n + 1$ different sequences of elementary expansions followed by an elementary collapse.  One such sequence is illustrated in the following picture

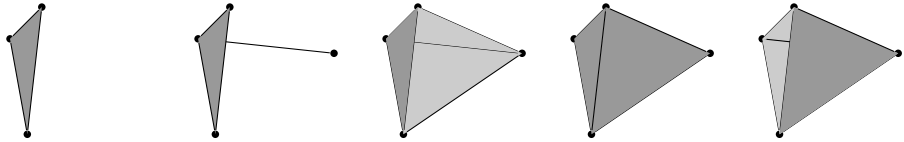
\begin{center}
  \begin{tikzpicture}
    \tikzstyle{vertex} = [circle, fill=black, inner sep=1pt]
    
    
    \node[vertex] (B) at (-1, 0, -1) {};
    \node[vertex] (C) at (0, 1, .5) {};
    \node[vertex] (D) at (0, -.5, 1) {};
    \fill[gray!80] (C.center) -- (D.center) -- (B.center) -- (C.center);

    \draw (C) -- (B) -- (D) -- (C);

    \node[vertex] (A) at (3.5, 0, -.5) {};
    \node[vertex] (B) at (1.5, 0, -1) {};
    \node[vertex] (C) at (2.5, 1, .5) {};
    \node[vertex] (D) at (2.5, -.5, 1) {};
    \draw (A) -- (B);
    \fill[gray!80] (C.center) -- (D.center) -- (B.center) -- (C.center);

    \draw (C) -- (B) -- (D) -- (C);

    \node[vertex] (A) at (6, 0, -.5) {};
    \node[vertex] (B) at (4, 0, -1) {};
    \node[vertex] (C) at (5, 1, .5) {};
    \node[vertex] (D) at (5, -.5, 1) {};
    \fill[gray!40] (A.center) -- (D.center) -- (B.center) -- (A.center);

    \draw (A) -- (B) -- (D) -- (A);
    \draw (B) -- (C);
    \draw (A) -- (C);
    \fill[gray!40] (A.center) -- (B.center) -- (C.center) -- (A.center);
    \draw (D) -- (C);
    \fill[gray!80] (B.center) -- (D.center) -- (C.center) -- (B.center);

    \node[vertex] (A) at (8.5, 0, -.5) {};
    \node[vertex] (B) at (6.5, 0, -1) {};
    \node[vertex] (C) at (7.5, 1, .5) {};
    \node[vertex] (D) at (7.5, -.5, 1) {};

    \draw (A) -- (B) -- (D) -- (A);
    \draw (B) -- (C);
    \fill[gray!80] (B.center) -- (D.center) -- (C.center) -- (B.center);
    \fill[gray!80] (A.center) -- (D.center) -- (C.center) -- (A.center);
    \draw (D) -- (C);
    \draw (A) -- (C);

    \node[vertex] (A) at (11, 0, -.5) {};
    \node[vertex] (B) at (9, 0, -1) {};
    \node[vertex] (C) at (10, 1, .5) {};
    \node[vertex] (D) at (10, -.5, 1) {};

    \draw (A) -- (B) -- (D) -- (A);
    \draw (B) -- (C);
    \fill[gray!40] (A.center) -- (D.center) -- (B.center) -- (A.center);
    \fill[gray!40] (A.center) -- (B.center) -- (C.center) -- (A.center);
    \draw (D) -- (C);
    \draw (A) -- (B);
    \fill[gray!80] (A.center) -- (D.center) -- (C.center) -- (A.center);

    \draw (A) -- (C);

  \end{tikzpicture}
\end{center}

Once again, we begin by embedding $\Delta^n$ in $\Delta^{n + 1}$ by identifying $\Delta^{n + 1}$ with $[e_{\star}, e_0, \dots, e_n]$.  We shall again use primes to indicate the face $\sigma' \subset \sigma$ in each elementary expansion.  The sequence of elementary expansion begins by choosing a vertex $e_j$ from $e_0, \dots, e_n$ and adding the edge $[e_{\star}', e_j]$.  Then for each vertices $e_k$ in $e_0, \dots, e_{j - 1}, e_{j + 1}, \dots, e_n$ we add the $2$-simplex $[e_{\star}', e_j, e_k']$.  Then for each pair of vertices $e_{k_1}, e_{k_2}$ in $e_0, \dots, e_{j - 1}, e_{j + 1}, \dots, e_n$, we add the $3$-simplex $[e_{\star}', e_j, e_{k_1}', e_{k_2}']$.  We continue inductively and end up adding an $l$-simplex for each of the $\binom{n}{l - 1}$ choices of $l - 1$ vertices.

Every simplex in the simplicial complex $\Delta^{n + 1}$ will be added by this procedure.  Clearly, any simplex containing both $e_{\star}$ and $e_j$ will be added.  Any simplex not containing $e_{\star}$ is already present from $\Delta^n$.  A simplex $[e_{\star}, e_{k_1}, \dots, e_{k_p}]$, which does not contain $e_j$, is added with $[e_0, e_j, e_{k_1}, \dots, e_{k_p}]$ in elementary expansion.

For each elementary expansion the inclusion map is given by
\begin{align*}
  i^{[e'_{\star}, e_j, e'_{k_1}, \dots, e'_{k_p}]}(\widehat{\tau}) =
  \begin{cases}
    \reallywidehat{[e_{\star}, e_j, e_{k_1}, \dots, \widehat{e_{k_l}}, \dots, e_{k_p}]} \\
    + (-1)^{l + 1}\reallywidehat{[e_\star, e_{k_1}, \dots, e_{k_p}]} & \text{if $\tau = [e_{\star}, e_j, e_{k_1}, \dots, \widehat{e_{k_l}}, \dots, e_{k_p}]$}\\
    \reallywidehat{[e_j, e_{k_1}, \dots, e_{k_p}]} \\
    + \reallywidehat{[e_{\star}, e_{k_1}, \dots, e_{k_p}]} & \text{if $\tau = [e_j, e_{k_1}, \dots, e_{k_p}]$}\\
    \widehat{\tau} & \text{otherwise}
  \end{cases}
\end{align*}
the projection map is given by
\begin{align*}
  p^{[e'_{\star}, e_j, e'_{k_1}, \dots, e'_{k_p}]}(\widehat{\tau}) =
  \begin{cases}
    0 & \text{if $\tau = [e_{\star}, e_j, e_{k_1}, \dots, e_{k_p}]$}\\
    0 & \text{if $\tau = [e_{\star}, e_{k_1}, \dots, e_{k_p}]$}\\
    \widehat{\tau} & \text{otherwise}
  \end{cases}
\end{align*}
and the homotopy is given by
\begin{align*}
  a^{[e_{\star}', e_j, e_{k_1}', \dots, e_{k_p}']}(\widehat{\tau}) =
  \begin{cases}
    -\reallywidehat{[e_{\star}, e_{k_1}, \dots, e_{k_p}]} & \text{if $\tau = [e_{\star}, e_j, e_{k_1}, \dots, e_{k_p}]$}\\
    0 & \text{otherwise}
  \end{cases}
\end{align*}

The inclusion map for the sequence of elementary expansions is thus
\begin{align*}
  i_j^{\downarrow}(\widehat{\tau}) =
  \begin{cases}
    \reallywidehat{[e_j, e_{k_1}, \dots, e_{k_p}]} \\
    +\reallywidehat{[e_{\star}, e_{k_1}, \dots, e_{k_p}]} & \text{if $\tau = [e_j, e_{k_1}, \dots, e_{k_p}]$ for $p = 0, \dots, n$}\\
    \widehat{\tau} & \text{otherwise}
  \end{cases}
\end{align*}
the projection map is thus
\begin{align*}
  p_j^{\downarrow}(\widehat{\tau}) =
  \begin{cases}
    0 & \text{if $\tau = [e_{\star}, e_j, e_{k_1}, \dots, e_{k_p}]$ for $p = 0, \dots, n - 1$}\\
    0 & \text{if $\tau = [e_{\star}, e_{k_1}, \dots, e_{k_p}]$ for $p = 0, \dots, n$}\\
    \widehat{\tau} & \text{otherwise}
  \end{cases}
\end{align*}
and the homotopy is thus
\begin{align*}
  a_j^{\downarrow}(\widehat{\tau}) =
    \begin{cases}
    -\reallywidehat{[e_{\star}, e_{k_1}, \dots, e_{k_p}]} & \text{if $\tau = [e_{\star}, e_j, e_{k_1}, \dots, e_{k_p}]$}\\
    0 & \text{otherwise}
  \end{cases}
\end{align*}
Composing with the elementary collapse $[e_{\star}, e_0', \dots, e_n']$ gives
\begin{align*}
  \myoset{\star}{i}_j(\widehat{\tau}) &= p^{[e_{\star}, e_0', \dots, e_n']} i_j^{\downarrow}(\widehat{\tau}) \\
  &= 
  \begin{cases}
     (-1)^j\reallywidehat{[e_{\star}, e_0, \dots, \widehat{e_j}, \dots, e_n]} & \text{if $\tau = [e_0, \dots, e_n]$}\\ 
    \reallywidehat{[e_j, e_{k_1}, \dots, e_{k_p}]} \\
    + \reallywidehat{[e_{\star}, e_{k_1}, \dots, e_{k_p}]} & \text{if $\tau = [e_j, e_{k_1}, \dots, e_{k_p}]$ for $p = 0, \dots, n - 1$}\\
    \widehat{\tau} & \text{otherwise}
  \end{cases}
\end{align*}
and
because
\begin{align*}
  i^{[e_{\star}, e'_0, \dots, e'_n]}(\widehat{\tau}) =
  \begin{cases}
    \reallywidehat{[e_{\star}, e_0, \dots, \widehat{e_l}, \dots, e_n]} \\
    - (-1)^l\reallywidehat{[e_0, \dots, e_n]} & \text{if $\tau = [e_{\star}, e_0, \dots, \widehat{e_l}, \dots,  e_n]$}\\
    \widehat{\tau} & \text{otherwise}
  \end{cases}
\end{align*}
we have
\begin{align*}
  \myoset{\star}{p}_j(\widehat{\tau}) &= p_j^{\downarrow} \, i^{[e_{\star}, e'_0, \dots, e'_n]}(\widehat{\tau}) \\
  &= 
  \begin{cases}
    (-1)^l\reallywidehat{[e_0, \dots, e_n]} & \text{if $\tau = [e_{\star}, e_0, \dots, \widehat{e_l}, \dots,  e_n]$ for $l = 0, \dots, n$}\\
    \widehat{\tau} & \text{otherwise}
  \end{cases}
\end{align*}
and lastly
\begin{align*}
  \myoset{\star}{a}_j(\widehat{\tau}) &= p^{[e_{\star}, e'_0, \dots, e'_n]}h_j^{\downarrow}i^{[e_{\star}, e'_0, \dots, e'_n]}(\widehat{\tau})\\
  &= 
    \begin{cases}
    -\reallywidehat{[e_{\star}, e_{k_1}, \dots, e_{k_p}]} & \text{if $\tau = [e_{\star}, e_j, e_{k_1}, \dots, e_{k_p}]$ for $p = 0, \dots, n - 1$}\\
    0 & \text{otherwise}
  \end{cases}
\end{align*}

Let $\myoset{\star}{\iota} = \frac{1}{n + 1}\sum_j\myoset{\star}{\iota}_j$, let  $\myoset{\star}{p} = \myoset{\star}{p}_j$ and let $\myoset{\star}{a} = \frac{1}{n + 1}\sum_j \myoset{\star}{a}_j$.  Because we have composed a zigzag of deformation retractions, we have $d\myoset{\star}{a} + \myoset{\star}{a}d = 1 - \myoset{\star}{\iota}\myoset{\star}{p}$.  Remarkably it is still the case the $\myoset{\star}{p}\myoset{\star}{\iota} = 1$.  Furthermore,
\begin{mythm}
  The deformation retraction
  \begin{align*}
    \myDR{C^{\bullet}(\Delta^n)}{C^{\bullet}(\star \Delta^n)}{\myoset{\star}{\iota}}{\myoset{\star}{p}}{\myoset{\star}{a}}
  \end{align*}
  is equal to the stellar subdivision on cochains deformation retraction of Theorem \ref{thm:stell-subd-cochains}.  That is $\myoset{\star}{\iota} = i^{\star}$, $\myoset{\star}{p} = p^{\star}$ and $\myoset{\star}{a} = a^{\star}$.
\end{mythm}

For simplicity, we have chosen in this section to focus on the construction for the stellar subdivision $\star \Delta^n = \star_{1, \dots, n} \Delta^n$.

We indicate how one would proceed in constructing the more general stellar subdivision $\star_{i_0, \dots, i_k} \Delta^n$.  For each $i_{\alpha} \in \{i_0, \dots, i_k\}$, we follow the same sequence of elementary expansions followed by an elementary collapse as specified above.  Each of these $k + 1$ sequences constructs $\star \Delta^n$ from $\Delta^n$.  We succeed each sequence with addition elementary collapses of all the simplices $[e_{\star}, e_{j_0}', \dots, e_{j_l}']$ containing $[e_{\star}, e_{i_0}, \dots, e_{i_k}]$ starting in dimension $n - 1$ and moving to lower dimensions.  Let $\myoset{\star}{\iota}_{\alpha}$, $\myoset{\star}{p}_{\alpha}$, and $\myoset{\star}{a}_{\alpha}$ be the inclusion, projection, and homotopy respectively that results from composing the specified sequence of elementary expansions followed by elementary collapses.  Let $\myoset{\star}{\iota} = \frac{1}{k + 1} \sum_{\alpha} \myoset{\star}{\iota}_{\alpha}$, let $\myoset{\star}{p} = \myoset{\star}{p}_{\alpha}$, and let $\myoset{\star}{a} = \frac{1}{k + 1} \sum_{\alpha} \myoset{\star}{a}_{\alpha}$.  In conclusion, we claim that these maps form a deformation retraction that is equal to the stellar subdivision deformation retraction that was constructed in Theorem \ref{thm:stell-subd-cochains}.

\section*{Acknowledgments}
Special thanks to Kevin Costello for helpful discussions.  This research was supported in part by Perimeter Institute for Theoretical Physics.  Research at Perimeter Institute is supported by the Government of Canada through the Department of Innovation, Science and Economic Development Canada and by the Province of Ontario through the Ministry of Economic Development, Job Creation and Trade.
\printbibliography
\end{document}